\documentclass[11pt,reqno]{amsproc}

\title[]{On the Nernst-Planck-Navier-Stokes system}

\author{Peter Constantin}
\address{Department of Mathematics, Princeton University, Princeton, NJ 08544}
\email{const@math.princeton.edu}

\author{Mihaela Ignatova}
\address{Department of Mathematics, Princeton University, Princeton, NJ 08544}
\email{ignatova@math.princeton.edu}

\usepackage[margin=1in]{geometry}
\usepackage{amsmath, amsthm, amssymb}
\usepackage{times}
\usepackage{color}
\usepackage{hyperref}

\newcommand{\pa}{\partial}
\newcommand{\la}{\label}
\newcommand{\fr}{\frac}
\newcommand{\na}{\nabla}
\newcommand{\be}{\begin{equation}}
\newcommand{\ee}{\end{equation}}
\newcommand{\ba}{\begin{array}{l}}
\newcommand{\ea}{\end{array}}
\newcommand{\Rr}{{\mathbb R}}
\renewcommand{\c}{{\widetilde{c}}}

\newcommand{\beg}{\begin}

\renewcommand{\div}{{\mbox{div}\,}}
\newcommand{\D}{\Delta}
\newcommand{\T}{T_K}
\date{today}
\begin{document}
\begin{abstract}
We consider ionic electrodiffusion in fluids, described by the Nernst-Planck-Navier-Stokes system in bounded domains, in two dimensions, with Dirichlet boundary conditions for the Navier-Stokes and Poisson equations, and blocking (vanishing normal flux) or selective (Dirichlet) boundary conditions for the ionic concentrations. We prove global existence and stability results for large data. 
\end{abstract}
\keywords{electroconvection, ionic electrodiffusion, Poisson-Boltzmann, Nernst-Planck, Navier-Stokes, relative entropy, Kullback-Leibler}

\noindent\thanks{\em{ MSC Classification:  35Q30, 35Q35, 35Q92.}}

\maketitle

\section{Introduction}
We consider electrodiffusion of ions in fluids in the presence of boundaries.
Ions of different valences carry charges  and diffuse under the influence of an electric potential, their own concentration gradients and a fluid flow. The fluid is forced by the electric forces created by the ions. The situation is described by the Nernst-Planck equations 
\be
\pa_t c_{i} + \div j_{i}\la {np} =0
\ee
where $c_i$ are the $i$-th ionic species concentrations, $i=1,\dots N$,
and where the fluxes $j_i$ are given by
\be
j_{i} = uc_{i} - D_{i}\na c_{i} - D_{i} \fr{ez_{i}}{k_BT}c_{i}\na \Psi.
\la{jpm}
\ee
The ion concentrations $c_{i} = c_{i}(x,t)$ are nonnegative functions, with $x$ representing position, $x\in \Omega\subset\Rr^d$, an open bounded set with smooth, orientable boundary, and $t$ representing time, $t\ge 0$. The domain is connected but not necessarily simply connected. The velocity $u=u(x,t)$ is a divergence-free field. $D_{i}$ are positive constant diffusivities ($D_i>0$, possibly different one from the other), $e$ is elementary charge, $z_i$ are valences ($z_i\in\Rr$, unrestricted), $k_B$ is Boltzmann's constant and $T$ is temperature. The potential $\Psi$ solves a Poisson equation
\be
-\varepsilon \D \Psi = {\widetilde{\rho}}
\la{poi}
\ee
in $\Omega$. The function ${\widetilde{\rho}}$ is the charge density,
\be
{\widetilde{\rho}} = e\sum_{i=1}^N z_ic_i 
\la{wrho}
\ee
and $\varepsilon$ is a positive constant, the dielectric permittivity of the solvent. The velocity $u$ obeys the Navier-Stokes equations 
\be
\pa_t u + u\cdot\na u -\nu\Delta u + \na p = {\widetilde{\rho}} {\widetilde{E}}
\la{nseqns}
\ee
in $\Omega$ with the divergence-free condition
\be
\na\cdot u = 0
\la{divu}
\ee
and with ${\widetilde{E}}$ the electric field
\be
{\widetilde{E}} = -\na \Psi.
\la{wele}
\ee
Here $\nu>0$ is the kinematic viscosity and $p$ the pressure. There are two kinds of boundary conditions for the ionic concentrations. The vanishing of all normal fluxes
\be
(j_{i}\cdot n)_{\left |\right. \pa\Omega} = 0, \quad i=1,\dots N,
\la{jbc}
\ee
where $n$ is outer normal at the boundary of $\Omega$, is termed ``blocking boundary conditions''. These boundary conditions model situations in which boundaries are impermeable: the ions are not allowed to cross them. 

Different boundary conditions are termed ``selective'' or ``permselective''. They model situations in which some ionic species are selectively crossing some boundaries, while being blocked from crossing others. In this case  $M\le N$ of the ionic concentrations have mixed Dirichlet - no-flux boundary conditions, and the rest of the ionic species ($i= M+1, \dots, N$) have blocking boundary conditions (\ref{jbc}).
\be
\left\{
\ba
{c_i}_{\left |\right. S_i} = \gamma_i, \quad (j_{i}\cdot n)_{\left |\right. \pa\Omega\setminus S_i} = 0, \quad i= 1, \dots M,\\
(j_{i}\cdot n)_{\left |\right. \pa\Omega} = 0, \quad i= M+1,, \dots, N.
\ea
\right.
\la{gammaibc}
\ee
where $S_i\subset \pa\Omega$ are portions of the boundary for $i=1, \dots, M$, and $\gamma_i>0$ are positive constants.  The subsets $S_i$ can be quite general: they do not need to be connected, nor do they need to be distinct from one another as $i$ varies. 

The electric potential satisfies Dirichlet boundary conditions
\be
\Psi_{\left | \right. \pa \Omega} = V
\la{dir}
\ee
where $V(x)$ are imposed voltages (the boundary $\pa \Omega$ need not be connected).
We normalize the potential by introducing $\Phi$,
\be
\Phi = \fr{e}{k_BT}\Psi,
\la{phi}
\ee 
(we depart somewhat from customary normalizations which include a valence), and denoting
\be
\rho = \sum_{i=1}^Nz_ic_i
\la{rho}
\ee
the NPNS system is therefore
\be
(\pa_t + u\cdot\na)c_{i} = D_{i}\div(\na c_{i} + z_i c_{i}\na\Phi) =
D_{i}\div(c_{i}\na(\log c_{i} + z_{i}\Phi))
\la{cpmeq}
\ee
together with
\be 
-\epsilon\Delta \Phi = \rho
\la{poiphi}
\ee
and the forced Navier-Stokes equations
\be
\pa_t u + u\cdot\na u -\nu \D u + \na p = -k_BT\rho\na\Phi, 
\la{nse}
\ee
\be
\na\cdot u = 0,
\la{divuz}
\ee
with
\be
\epsilon = \fr{\varepsilon k_B T}{e^2}.
\la{eps}
\ee
We did not rescale the equations, we just slightly changed the dependent variables potential and charge density. We note that $\epsilon$ is essentially a length squared,
\be
\epsilon = c_0(\sum_{i=1}^Nz_i^2)\lambda_D^2
\la{deb}
\ee
where $\lambda_D$ is the Debye screening length and $c_0$ a reference bulk concentration of ions.  The boundary conditions for $u$ are homogeneous Dirichlet,
\be
u_{\left |\right. \pa \Omega} = 0,
\la{ubc}
\ee
and the blocking boundary conditions (\ref{jbc}) for $c_i$ thus become
\be
c_{i}\pa_n(\log c_{i} + z_{i}\Phi)_{\left |\right. \pa\Omega} = 0,
\la{cbc}
\ee
where
\be
\pa_n = n\cdot\na
\la{pan}
\ee
is normal derivative at the boundary.
The boundary condition for $\Phi$ is
\be
\Phi_{\left |\right. \pa\Omega} = W = \frac{e}{k_BT}V,
\la{phibc}
\ee
with $W=W(x)$ a given smooth enough function of space. We distinguish between two kinds of selective boundary conditions for the concentrations $c_i$. The first, which we term ``uniform selective'', require not only the $\gamma_i$ to be constant (in space and time) but also that the boundary voltage $W(x)$ to be constant on the portions $S_i$ of the boundary where $\gamma_i$ are prescribed,
\be
W(x)_{\left |\right.\; S_i} = w_i.
\la{Wwi}
\ee 
The rest of selective boundary conditions we term  ``general selective''. In their case $W(x)$ may be an arbitrary (smooth enough) function of space.

The Boltzmann steady states are defined to be
\be
c_{i}^*(x) = \frac{e^{-z_i \Phi^*(x)}}{Z_{i}}
\la{cstar}
\ee
with $Z_{i}>0$ constants (which may depend on $\Phi^*$). The function $\Phi^*(x)$ is time independent and obeys the semilinear elliptic equation 
\be 
-\epsilon \D \Phi^* = \rho^*
\la{poistar}
\ee
with
\be
\rho^* = \sum_{i=1}^Nz_ic_i^*
\la{rhostar}
\ee
and with boundary condition
\be
\Phi^*_{\left |\right. \pa\Omega} = W.
\la{phistarbc}
\ee
This equation is known as the Poisson-Boltzmann equation. Let us observe that $c_{i}^*$, $\Phi^*$ are steady solutions of the NPNS system with $u=0$. Indeed, in this situation the forcing term in the Navier-Stokes equations (\ref{nse}) is a gradient and it can be included in the pressure, while the time independent equations (\ref{cpmeq}) are satisfied.

The NPNS system is nonlinear, and the blocking boundary conditions are nonlinear and nonlocal. The physical and biophysical applications of the system are extremely broad, and the system has been investigated extensively in the physical literature. An introduction to some of the basic physical and mathematical issues can be found in \cite{rubibook}. While blocking boundary conditions lead to stable configurations, instabilities occur for selective boundary conditions. These have been studied in simplified models mathematically and numerically (\cite{rubizaltz}, \cite{zaltzrubi}) and observed in physical experiments \cite{rubinstein}. A recent numerical study, which partly motivated ours, \cite{davidson}, discussed additional ``patterned'' boundary conditions, and described the effect of the geometry of nonuniform boundary  conditions on the instabilities.
The numerical study is performed in a strip, with periodic lateral boundary conditions. There are two ionic species, anions and cations, and the boundary conditions for anions are blocking while the boundary conditions for cations are selective. The boundary conditions for the electric potential are Dirichlet: a constant voltage is applied at one of the boundaries. The case when both boundaries for cations are selective  corresponds in our language to general selective boundary conditions: $N=2$,  $S_1 =\pa\Omega$ is formed by both the upper and the lower boundary, $c_1$ is constant on $S_1$, but $W$ is not, taking two different values. An interesting other case is one in which the upper boundary for cations is selective and the lower boundary is patterned with alternating segments of permeable and impermeable membranes. Both situations lead to instability and chaotic behavior, and correspond in our language to general selective boundary conditions. Interestingly, if the upper boundary is blocking, but the lower one is selective, or even patterned selective, then we are in situations which we call ``uniform'' selective, because the voltage is constant on the selective part of the boundary. These, and more complicated cases with many boundary components and many ion species are proved in this paper to be nevertheless unconditionally globally uniformly stable situations.

The mathematical study of the relevant semilinear elliptic equations is classical (\cite{keller}, \cite{friedman}). The coupled NPNS system is semilinear parabolic, so its local well posedness is not unexpected. The issue is whether or not solutions exist globally and what is their asymptotic behavior. This issue is mostly a question of boundary conditions, although dimensionality enters as well. Global existence and stability of solutions of the Nernst-Planck equations, uncoupled to fluids has been obtained in several situations in \cite{biler}, \cite{choi}, \cite{gajewski} for blocking boundary conditions. Local existence for the system coupled to the Navier-Stokes equations in the whole space was obtained in \cite{jerome} and global existence of weak solutions in 3D with blocking boundary conditions was obtained in \cite{jeromesacco} and in \cite{fischer}. The global existence and stability of the system in 2D has been studied in \cite{bothe} with blocking boundary conditions for the ions and a Robin boundary condition for the electric potential. The method of proof and the result of \cite{bothe} do not apply to the case of Dirichlet boundary conditions for the potential. Global existence for small data and forces was obtained in \cite{ryham} and \cite{schmuck}.

In this paper we prove global existence for both blocking and selective boundary conditions for the ionic concentrations, in two spatial dimensions, for arbitrary data. In the cases of blocking boundary conditions and in the case of uniform selective boundary conditions we prove unconditional global stability: for arbitrary large initial data, valences, voltages, species diffusivities, dielectric constant and arbitrary Reynolds numbers, the solutions converge as time tends to infinity to unique selected Boltzmann states. The Boltzmann states are uniquely determined by the initial average concentrations of the species and boundary conditions. The Navier-Stokes equations are forced, and the forces converge in time to potential forces, but they are not, in general, potential forces at any finite time. Thus the fact that the attractor is a singleton (per leaf) is nontrivial, and it follows from the remarkable structure of the equations: The sum of natural relative entropies (or Kullback-Leibler divergences), relative to Boltzmann states, together with the mean-square gradient difference of electrical potential and the kinetic energy of the fluid decays in time.

The paper is organized as follows: in Section 2 we describe the dissipative structure, in Section 3 we give a priori bounds and decay to Boltzmann states for blocking boundary conditions, in Section 4 we describe the stability of uniform selective boundary conditions and in Section 5 we describe the global existence for the general selective boundary conditions. Appendix A is devoted to the Poisson-Boltzmann equations, and Appendix B to a proof of local existence.

\section{Dissipative Structure}
Let us define the energy 
\be
{\mathcal E} = {\mathcal E}(c_i, \Phi; c_i^*, \Phi^*) = \int_{\Omega}\left[\sum_{i=1}^N E_i c_i^*  +\fr{1}{2}(\rho-\rho^*)(\Phi-\Phi^*)\right]dx.
\la{enew}
\ee
This energy is relative to some fixed selected Boltzmann states, 
\be
c_i^*(x) = Z_{i}^{-1}e^{-z_i\Phi^*(x)}
\la{cistar}
\ee
obeying the Poisson-Boltzmann equation (\ref{poistar}) with boundary conditions (\ref{phistarbc}) and charge density
\be
\rho^*(x) = \sum_{i=1}^N z_ic_i^*(x).
\la{rostarn}
\ee
Above we used $\rho =\sum z_ic_i$ (see \ref{rho}), and denoted $E_i$ by
\be
E_i = \fr{c_i}{c_i^*}\log\left(\fr{c_i}{c_i^*}\right) - \fr{c_i}{c_i^*} + 1.
\la{ei}
\ee
We have the relations
\be
\fr{\pa (E_i c_i^*)}{\pa c_i} = \log\left(\fr{c_i}{c_i^*}\right)
\la{paeici}
\ee
and
\be
\fr{\pa (E_i c_i^*)}{\pa c_i^*} = 1 - \fr{c_i}{c_i^*}.
\la{paeicistar}
\ee
The potential $\Phi$ in ${\mathcal E}$ is computed solving the Poisson problem (\ref{poiphi})
\be
-\epsilon\Delta \Phi = \rho
\la{poin}
\ee
with boundary condition (\ref{phibc})
\be
\Phi_{\left | \right. \;\pa\Omega} = W.
\la{poinbc}
\ee
Computing the first variations (Fr\'{e}chet derivatives) of ${\mathcal E}$ gives the densities
\be
\fr{\delta \mathcal E}{\delta c_i} = \log\left(\fr{c_i}{c_i^*}\right) + z_i(\Phi-\Phi^*)
\la{deleici}
\ee
because
\be
\fr{1}{2}(\rho - \rho^*)(\Phi-\Phi^*) = \fr{1}{2\epsilon}(\rho-\rho^*)(-\D_D)^{-1}(\rho-\rho^*)
\la{potentialene}
\ee
where $(-\D_D)^{-1}$ is the inverse Laplacian with homogeneous Dirichlet boundary conditions, which is a selfadjoint operator, and because
\be
\fr{\pa \rho}{\pa c_{i}} =  z_{i}.
\la{rhoc}
\ee
Note that in view of (\ref{cistar}) we have
\be
\fr{\delta \mathcal E}{\delta c_i} = \log c_i + z_i\Phi + \log Z_i
\la{deleicin}
\ee
and therefore the equations (\ref{cpmeq})
\be
D_t c_i = D_i\div\left(c_i\na\left(\log c_i + z_i\Phi\right)\right)
\la{dtci}
\ee
are, in view of the relation (\ref{deleicin}),  the same as
\be
D_t c_i = D_i\div\left(c_i\na\left(\fr{\delta\mathcal E}{\delta c_i}\right)\right).
\la{cieq}
\ee
We denoted above by $D_t$ the material derivative 
\be
D_t =\pa_t + u\cdot\na
\la{dt}
\ee
with respect to the time dependent, divergence-free velocity $u$. 
This is a fundamental property of the Nernst-Poisson system. The variational structure is obtained using only the fact that $\Phi$ and $\Phi^*$ obey the same boundary conditions. 
Defining the energy density by 
\be
E = \sum_{i=1}^N E_ic_i^* + \fr{1}{2}(\rho-\rho^*)(\Phi-\Phi^*)
\la{densen}
\ee
we compute $D_t (E_ic_i^*)$ using (\ref{paeici}), (\ref{paeicistar}) and (\ref{cistar}): 
\be
\ba
D_t (E_ic_i^*) = \log\left(\fr{c_i}{c_i^*}\right)D_tc_i + D_tc_i^* - c_i D_t\log c_i^*\\
= \log\left(\fr{c_i}{c_i^*}\right)D_tc_i + D_tc_i^* + z_ic_iD_t \Phi^*.
\ea
\la{dteicistar}
\ee
Adding we obtain
\be
D_t \left (\sum_{i=1}^N E_ic_i^*\right) = \sum_{i=1}^N \log\left(\fr{c_i}{c_i^*}\right)D_tc_i + \sum_{i=1}^N D_t c_i^* + \rho D_t\Phi^*.
\la{dtentr}
\ee
In view of (\ref{deleici}) we have thus
\be
D_t \left (\sum_{i=1}^N E_ic_i^*\right) = \sum_{i=1}^N \fr{\delta\mathcal E}{\delta c_i}D_tc_i + \sum_{i=1}^N D_tc_i^* - (\Phi-\Phi^*)D_t\rho + \rho D_t\Phi^*.
\la{dtentrop}
\ee
Therefore
\be
D_t E = \sum_{i=1}^N \fr{\delta\mathcal E}{\delta c_i}D_tc_i + \sum_{i=1}^N D_tc_i^* + P
\la{dteP}
\ee
where
\be
P = \fr{1}{2}D_t[(\rho - \rho^*)(\Phi-\Phi^*)] + \rho D_t\Phi^* - (\Phi-\Phi^*)D_t\rho,
\la{P}
\ee 
and thus
\be
P = \fr{1}{2}D_t\left[\rho \Phi^* -\rho^*\Phi + \rho^*\Phi^*\right] + 
\fr{1}{2}\rho D_t\Phi - \fr{1}{2}(D_t\rho)\Phi.
\la{pcalc}
\ee
Now we claim that
\be
P = \rho u\cdot\na \Phi + Q
\la{Pq}
\ee
and
\be
\int_{\Omega} Q dx = 0
\la{intq}
\ee
for all $t$. Indeed,
\be
\ba
Q =\fr{1}{2}D_t\left[\rho \Phi^* -\rho^*\Phi + \rho^*\Phi^*\right] -\fr{1}{2}\div(u\rho\Phi) + \fr{1}{2}(\rho \pa_t\Phi -\Phi\pa_t\rho)\\
=\fr{1}{2}[\Phi^*\pa_t\rho -\rho^*\pa_t\Phi + \rho \pa_t\Phi -\Phi\pa_t\rho] +
\fr{1}{2}\div[u\left(\rho \Phi^* -\rho^*\Phi + \rho^*\Phi^* -\rho\Phi\right)]\\
= \fr{1}{2}[(\rho-\rho^*)\pa_t\Phi - (\Phi-\Phi^*)\pa_t\rho] + \fr{1}{2}\div[u\left(\rho \Phi^* -\rho^*\Phi + \rho^*\Phi^* -\rho\Phi\right)]\\
=\fr{1}{2}[(\rho-\rho^*)\pa_t(\Phi-\Phi^*) - (\Phi-\Phi^*)\pa_t(\rho-\rho^*)]  + \fr{1}{2}\div\left[u(\rho +\rho^*)( \Phi^* -\Phi)\right] 
\ea
\la{qint}
\ee
where we used that 
\be
\pa_t\Phi^* = \pa_t \rho^* = 0.
\la{startime}
\ee
Thus
\be
Q = \fr{1}{2\epsilon}[(\rho-\rho^*)(-\D_D)^{-1}(\pa_t(\rho-\rho^*)) -
((-\D_D)^{-1}(\rho -\rho^*))(\pa_t(\rho-\rho^*))] + \fr{1}{2}\div (uq)
\la{Qq}
\ee
with 
\be
q = (\rho +\rho^*)( \Phi^* -\Phi).
\la{q}
\ee
The fact that (\ref{intq}) holds follows from the facts that $(-\D_D)^{-1}$ is
selfadjoint and the fact that $u$ is divergence-free and has vanishing normal component on the boundary of $\Omega$. No boundary conditions on $c_i$ are used.   
We have thus

\be
D_tE =   \sum_{i=1}^N \fr{\delta\mathcal E}{\delta c_i}D_tc_i + \sum_{i=1}^N D_tc_i^* + \rho u\cdot\na\Phi + Q
\la{dteq}
\ee
where $Q$ satisfies (\ref{intq}). Consequently, we have
\be
D_t E = \sum_{i=1}^N \fr{\delta\mathcal E}{\delta c_i}D_tc_i  - F\cdot u + R
\la{dter}
\ee
where
\be
F = -\rho \na\Phi
\la{F}
\ee
and with
\be
R =  \sum_{i=1}^N D_tc_i^* + Q.
\la{R}
\ee
In view of (\ref{intq}) and of
\be
\pa_t c_{i}^* = 0
\la{dtcstars}
\ee
we have that $R$ satisfies 
\be
\int_{\Omega}R(x,t)dx = 0
\la{intr}
\ee
for all $t$. We stress that no boundary conditions for $c_i$ were used so far. Now we use the coupling to the 2D Navier-Stokes equations whose kinetic energy is forced by $F$. Adding the energy balance in the Navier-Stokes equations multiplied by $\fr{1}{k_BT}$ we obtain from (\ref{cieq}), (\ref{dter}) and (\ref{intr}) after integration by parts
\be
\la{endissgenbc}
\fr{d}{dt}\left [\fr{1}{2k_BT}\int_{\Omega}|u|^2 dx + \mathcal E\right] = -{\mathcal D}  -\fr{\nu}{k_B T}\int_{\Omega}|\na u|^2dx  + \sum_{i=1}^ND_i\int_{\pa\Omega}c_i\fr{\delta\mathcal E}{\delta c_i}\pa_n\left(\fr{\delta\mathcal E}{\delta c_i}\right)dS
\ee
where
\be
{\mathcal D} = \sum_{i=1}^N D_i\int_{\Omega}c_i\left|\na \fr{\delta \mathcal E}{\delta c_i}\right|^2 dx.
\la{mathcaldn}
\ee
If  blocking boundary conditions (\ref{jbc}) are employed, then
\be
\pa_n\left(\fr{\delta\mathcal E}{\delta c_i}\right)_{\left |\right. \pa\Omega} = 0
\la{blockdele}
\ee
no matter what Boltzmann states are considered, in view of (\ref{deleicin}) and
the fact that $Z_i$ are constant in space (and time, of course). We recall that in this case $W$ is an arbitrary (smooth enough) function.

If Dirichlet boundary conditions (\ref{gammaibc}) are used in the case of uniform selective boundary conditions, then we choose
\be
Z_i = \left(\gamma_i e^{z_i w_i}\right)^{-1}, \quad {\mbox{for}}\; i=1,\dots M
\la{zichoice}
\ee
where we recall that
\be
w_i = W_{\left|\right. \; S_i}
\la{wi}
\ee
are assumed to be constant on $S_i$.  The rest of $Z_i$, $i = M+1,\dots, N$ are arbitrary and $W$ may vary in space on the rest of the boundary $\pa\Omega \setminus\cup_{i=1}^MS_i$. In this case we have, in view of  (\ref{deleicin}) and (\ref{gammaibc})
\be
\fr{\delta\mathcal E}{\delta c_i}\pa_n\left(\fr{\delta\mathcal E}{\delta c_i}\right)_{\left|\right .\; \pa\Omega} =0
\la{zerobc}
\ee
for all $i=1, \dots, N$.

\beg{thm}\la{diss}
Let $c_i>0$ solve the 2D Nernst-Planck-Navier-Stokes equations (\ref{rho}), (\ref{cpmeq}), (\ref{poiphi}), (\ref{nse}), (\ref{divuz}) with Dirchlet boundary conditions for the Navier-Stokes velocity (\ref{ubc}) and the electric potential (\ref{phibc}) and either blocking (vanishing normal flux) boundary conditions (\ref{jbc}) or uniform selective (constant Dirichlet and vanishing normal flux) boundary conditions (\ref{gammaibc})  for the ion concentrations.
Let ${\mathcal E}$ be defined in (\ref{enew}) with respect to arbitrary Boltzmann states in the case of blocking boundary conditions, and with respect to Boltzmann states selected by (\ref{zichoice}) for uniform selective boundary conditions. Then
\be
\fr{d}{dt}\left [\fr{1}{2k_BT}\int_{\Omega}|u|^2 dx + \mathcal E\right] = -{\mathcal D}  -\fr{\nu}{k_B T}\int_{\Omega}|\na u|^2dx
\la{gendiss}
\ee
holds for all $t>0$, where $\mathcal D$ is given by (\ref{mathcaldn}).
\end{thm}
\beg{rem}
The energy in the left hand side of (\ref{gendiss}) is non-negative. The energy is the sum of relative entropies (or Kullback-Leibler divergences) for the pairs $(c_i, c_i^*)$, the  square of the $H^{-1}$ norm of the difference of charge densities, and the kinetic energy of the fluid. It vanishes only if $c_i = c_i^*$, $\Phi=\Phi^*$ and $u=0$. The dissipation ${\mathcal D}$ also vanishes only at Boltzmann states. The dimension $d$ of space does not enter these calculations, and the only use of the Navier-Stokes equations is by considering $D_t$ as a derivation, and using the energy equality. In $d=3$, and for weak Leray solutions of the forced NSE, (\ref{gendiss}) holds with inequality rather than equality, for almost all time. The fact that (\ref{dter}) with (\ref{intr}) holds represents a mathematical confirmation that $F$ is the correct electrical forcing of Navier-Stokes or Stokes equations: no other force would have fulfilled its role. In other words, we could derive the form of $F$ by the requirement that Theorem {\ref{diss}} holds. The dimension of $k_BT $ is that of an energy, and (\ref{gendiss}) is dimensionally correct.
\end{rem}
\beg{rem} The right hand side of (\ref{gendiss}) is independent of the choice of reference Boltzmann state, in view of (\ref{deleicin}). This might seem puzzling, but is explained by the fact that the difference between two energies $\mathcal E_1$ and $\mathcal E_2$ corresponding to two different admissible choices of $Z_i$ is time independent. Indeed, this difference is the sum of time independent quantities and constant multiples of $\int_{\Omega} c_i(x,t)dx$ (for all $i$ in the case of blocking conditions and for $i= M+1, \dots, N$ for uniform selective boundary conditions) which are conserved under the evolution. This follows from the calculation below. 
Let $c_i^*, \Phi^*$ be the unique Boltzmann state corresponding to constants $Z_i$, and let $d_i^*, \Psi^*$ be the Boltzmann state corresponding to different constants $U_i>0$, which still satisfy the conditions (\ref{zichoice}) in the case of selective boundary conditions. Denote $q^* = \sum_{i=1}^Nz_id_i^*$. Let $E_1$ denote the energy density of the state $c_i, \Phi$, relative to the first Boltzmann state given by (\ref{densen}), and $E_2$ the energy density corresponding to the second state. The difference of densities is
\be
E_1-E_2 = \sum_{i=1}^Nc_i\log\left(\fr{d_i^*}{c_i^*}\right) + 
\fr{1}{2}\Phi\left(q^*-\rho^*\right) + \fr{1}{2}\rho\left(\Psi^*-\Phi^*\right) + \sum_{i=1}^N (c_i^*-d_i^*) + \fr{1}{2}\left( \rho^*\Phi^* - q^*\Psi^*\right).
\la{difeE}
\ee
Using (\ref{cistar}) and its analogue, we have
\be
\log\left(\fr{d_i^*}{c_i^*}\right) = z_i\left(\Phi^*-\Psi^*\right ) +  \log\left(\fr{Z_i}{U_i}\right),
\la{cstardstar}
\ee
and from (\ref{difeE}) it follows that
\be
E_1-E_2 = \sum_{i=1}^N c_i\log\left(\fr{Z_i}{U_i}\right) + \fr{1}{2}\rho\left(\Phi^*-\Psi^*\right) + \fr{1}{2}\Phi\left(q^*-\rho^*\right) +  \sum_{i=1}^N (c_i^*-d_i^*) + \fr{1}{2}\left( \rho^*\Phi^* - q^*\Psi^*\right).
\la{diffE}
\ee
Now we write 
\be
\Phi = \Phi_W + \Phi_0
\la{phis}
\ee
with
\be
-\epsilon\D\Phi_0 = \rho, \quad {\Phi_0}_{\left|\right . \; \pa\Omega} = 0,
\la{Phizero}
\ee
and
\be
-\epsilon \D\Phi_W = 0, \quad {\Phi_W}_{\left|\right . \; \pa\Omega} = W
\la{Phiw}
\ee
and rewrite (\ref{diffE}) as
\be
E_1-E_2 = \sum_{i=1}^N c_i\log\left(\fr{Z_i}{U_i}\right) + \fr{1}{2}\rho\left(\Phi^*-\Psi^*\right) + \fr{1}{2}\Phi_0\left(q^*-\rho^*\right) +  K^*
\la{diffeonetwo}
\ee
where
\be
K^* = \fr{1}{2}\Phi_W\left(q^*-\rho^*\right) +  \sum_{i=1}^N (c_i^*-d_i^*) + \fr{1}{2}\left( \rho^*\Phi^* - q^*\Psi^*\right)
\la{kstar}
\ee
is time independent. Therefore
\be
E_1-E_2 = \sum_{i=1}^N c_i\log\left(\fr{Z_i}{U_i}\right) + \fr{1}{2\epsilon}\rho(-\D_D)^{-1}\left(\rho^*-q^*\right) - \fr{1}{2\epsilon}\left((-\D_D)^{-1}\rho\right)\left(\rho^* - q^*\right) +  K^*,
\la{diffes}
\ee
and, integrating and using the selfadjointness of $(-\D_D)^{-1}$ we have
\be
\mathcal E_1- \mathcal E_2 = \sum_{i=1}^N\log\left(\fr{Z_i}{U_i}\right)\int_{\Omega}c_i(x,t)dx + \int_{\Omega}K^*dx.
\la{diffcales}
\ee
The integrals $\int_{\Omega} c_i(x,t)dx$ are time independent for all $i= 1, \dots, N$, if blocking conditions are used, and for $i= M+1, \dots , N$, if uniform selective boundary conditions are used. In the latter case, in view of (\ref{zichoice}), $\log\left(\fr{Z_i}{U_i}\right) = 0$ for $i=1,\dots, M$. Thus, it makes no difference which admissible Boltzmann state is chosen to define the energy for (\ref{gendiss}) to hold.
\end{rem}

\beg{rem} The decay of energy (\ref{gendiss}) implies that any time independent solution of the system is a Boltzmann state. The velocity vanishes and the gradient of pressure balances the electrical forces, which are a gradient, in steady state. It is interesting to note the fact that the electrical forces are a gradient only in steady state, not on the way to steady state. Because of this, the decay of velocity is non-trivial, as the Navier-Stokes are forced by a nondecaying force.
\end{rem}

\section{Global unconditional stability for blocking boundary conditions}
We consider the equations (\ref{dtci}) with boundary conditions 
\be
j_i\cdot n_{\left | \right .\; \pa\Omega} = 0
\la{jibc}
\ee
for the fluxes 
\be
j_i = uc_i -D_i (\na c_i  + z_ic_i \na\Phi).
\la{ji}
\ee
We first show that if $c_i(x,t)$ are positive at $t=0$, then they remain positive, as long as the solutions are regular. In order to show this we take a convex function $F:\Rr\to \Rr$ that is nonnegative, twice continuously differentiable, identically zero on the positive semiaxis, and strictly positive on the negative axis. We also assume
\be
F''(y)y^2 \le CF(y)
\la{fass}
\ee
with $C>0$ a fixed constant. Examples of such functions are
\be
F(y) = \left\{
\ba
y^{2m} \quad \quad {\mbox{for}}\quad y<0,\\
0 \quad \quad \quad {\mbox{for}}\quad y\ge 0
\ea
\right.
\la{Fm}
\ee
with $m>1$. (In fact $m=1$ works as well, although we have only $F\in W^{2,\infty}(\Rr)$ in that case.) We multiply the equation (\ref{dtci}) by $F'(c_i)$ and integrate by parts using (\ref{jibc}). We obtain
\be
\fr{d}{dt}\int_{\Omega}F(c_i)dx = -D_i\int_{\Omega} F''(c_i)\left[ |\na c_i|^2  + z_ic_i \na\Phi\cdot\na c_i\right]dx.
\la{intfeq}
\ee
Using a Schwartz inequality we have
\be
\fr{d}{dt}\int_{\Omega}F(c_i(x,t))dx \le \fr{CD_i}{2}z_i^2\|\na\Phi\|_{L^{\infty}(\Omega)}^2\int_{\Omega}F(c_i(x,t))dx.
\la{intfineq}
\ee
If $c_i(x,0)\ge 0$ then $F(c_i(x,0))=0$ and (\ref{intfineq}) above shows that 
$F(c_i(x,t))$ has vanishing integral. As $F$ is nonnegative, it follows that $F(c_i(x,t))= 0$ almost everywhere in $x$ and because $F$ does not vanish for negative values it follows that $c_i(x,t)$ is almost everywhere nonnegative.   
\beg{rem} The paper is arranged in what we consider to be a natural order, mostly based on the interest of the subject and on the flow of ideas, but not strictly on the logical order dictated by rigor. Thus, the positivity of solutions follows from approximations and local existence in which $\|\na\Phi\|_{L^2(L^{\infty})}$ is guaranteed to be finite and the initial integral of $F$ is finite. This remark applies, mutatis mutandis to the whole article. 
\end{rem}

We consider now a priori bounds on solutions. From (\ref{gendiss}) it follows that
\be
\sum_{i= 1}^N\int_{\Omega}\left(\fr{c_i(x,t)}{c_i^*(x)}\log\left(\fr{c_i(x,t)}{c_i^*(x)}\right) -\fr{c_i(x,t)}{c_i^*(x)} +1\right)c_i^*(x)dx \le {\mathcal E}(0) + \fr{1}{2k_B\T}\|u_0\|^2_{L^2(\Omega)}
\la{lloglb}
\ee
holds for all time $t$. Above and in what follows we will use $\T$ to denote temperature, which is a fixed constant, in order to avoid confusion with $T$ representing time. In view of (\ref{lloglb}) we know that
\be
\int_{\Omega} c_i(x,t) \log (c_i(x,t) + 2)dx \le C^*\left[{\mathcal E}(0) + \fr{1}{2k_B\T}\|u_0\|^2_{L^2(\Omega)}\right]
\la{cillogl}
\ee
holds for all $t$, with $C^*$ depending only on  bounds on $c_i^*$ and $z_i$.
Let us denote
\be
\Gamma = \left[{\mathcal E}(0) + \fr{1}{2k_B\T}\|u_0\|^2_{L^2(\Omega)}\right]
\la{Gamma}
\ee
a constant depending on the initial energy and velocity $L^2$ norm.   Consequently we have that
\be
\int_{\Omega}|\rho(x,t)|\log\left(|\rho(x,t)| + 2\right)dx \le C^*\Gamma
\la{rhob}
\ee
holds uniformly in time, with a slightly different $C^*$.
We consider the case $d=2$. 
\vspace{.5cm}

\noindent{\bf {Step 1}: $L^{\infty}$ bound on $\Phi$.} \\
From the Poisson equation 
\be
-\epsilon\D (\Phi-\Phi^*) = \rho-\rho^*
\la{poiphin}
\ee
with homogeneous Dirichlet boundary conditions we obtain that
\be
\|\Phi(\cdot, t)\|_{L^{\infty}(\Omega)} \le C^*\Gamma_1
\la{phinf}
\ee
uniformly in time, with $C^*$ depending on $\epsilon$, and the domain $\Omega$
and 
\be
\Gamma_1 = \Gamma  + \int_{\Omega}|\rho^*(x)|\log\left(|\rho^*(x)| + 2\right)dx.
\la{gammaone}
\ee
The proof of this fact follows from properties of the 
Green's function and from the fact that the Legendre transform of $x\log x-x+1$ defined on the semipositive axis is $e^x-1$, and consequently
\be
|\log|x-y| |\rho (y,t)|| \le |\rho(y,t)|\log|\rho(y,t)| - |\rho(y,t)| + e^{|\log|x-y||}.
\la{legendreineq}
\ee
This is an essential use of $d=2$. 

\vspace{.5cm}

\noindent {\bf Step 2: Local uniform $L^1(L^q)$ bounds for $c_i$.}\\
We exploit the fact that
\be
\int_0^T{\mathcal D}(t)dt<\infty.
\la{disint}
\ee
Because of (\ref{deleicin}) and (\ref{mathcaldn}) we have that 
\be
\int_0^T\int_{\Omega}c_i(x,t)\left|\na\log\left(c_i(x,t)e^{z_i\Phi(x,t)}\right)\right|^2dx dt \le {\mathcal E}(0) + \fr{1}{2k_B\T}\|u_0\|^2_{L^2(\Omega)} =\Gamma.
\la{dissphiexp}
\ee
Using the crucial information from the previous step that $\Phi$ is bounded a priori in $L^{\infty}$ (\ref{phinf}) we deduce that the useful auxiliary function
\be
\widetilde{c_i}(x,t) = c_i(x,t)e^{z_i\Phi(x,t)}
\la{tildeci}
\ee
obeys
\be
\int_0^T\int_{\Omega} {\widetilde{c_i}(x,t)}^{-1}\left |\na\widetilde{c_i}(x,t)\right|^2 dxdt \le C^*\Gamma e^{C^*\Gamma_1} = C^*\Gamma_2.
\la{natildeci}
\ee
Together with (\ref{phinf}) and (\ref{cillogl}), this implies that 
$\sqrt{\widetilde{c}_i}\in L^2(0,T; H^1(\Omega))$, and thus 
$\widetilde{c_i}\in L^1(0,T; L^q(\Omega))$ for any $q\in [1, \infty)$, with bounds depending only on the initial energy and growing linearly in $T$. More precisely, we have that
$\sqrt{\widetilde{c_i}}\in L^{\infty}(0,T; L^2(\Omega))$ and $\na\sqrt{\widetilde{c_i}}\in L^2(0,T; L^2(\Omega))$ and so, for any interval $[t_0,t_0+\tau]\subset [0,T]$ we have 
\be
\int_{t_0}^{t_0+\tau}\|\sqrt{\widetilde{c_i}(t)}\|_{H^1(\Omega)}^2dt \le C^*\Gamma_2(1+\tau)
\la{h1sqrt}
\ee
with $C^*$ independent of initial data and of time. Time enters in the right-hand side of the estimate because unlike its gradient which is mean square time integrable, the $\sqrt{\widetilde c_i}$ norms are bounded but not decaying in time.
Returning to $c_i$ and using again (\ref{phinf}) we obtain
\be
\int_{t_0}^{t_0+\tau}\|c_i(t)\|_{L^q(\Omega)}dt \le C^* \Gamma_3(1+\tau)
\la{cilq1}
\ee
with $\Gamma_3$ depending only on the initial energy and velocity $L^2$ norm via $\Gamma$,  and on $\rho^*$ via $\Gamma_1$. The constant $C^*$ depends on $q$ because we used embedding theorems. 

\vspace{.5cm}

\noindent{\bf{Step 3: Local uniform bounds for $c_i$ in $L^2(L^2)$.}}\\
In view of the fact that $\sqrt{\widetilde{c_i}}$ is bounded in $L^2(0,T; H^1(\Omega))$  (\ref{natildeci}) we can interpolate using Ladyzhenskaya (Gagliardo-Nirenberg) inequalities 
\be
\int_{\Omega}|\sqrt{\widetilde{c_i}(x,t)}|^4dx \le C\left(\int_{\Omega}|\sqrt{\widetilde{c_i}(x,t)}|^2dx\right)\|\sqrt{\widetilde{c_i}(t)}\|^2_{H^1(\Omega)},
\la{sqrtlady}
\ee
and, in view of the fact that $\sqrt{\widetilde{c_i}}\in L^{\infty}(0,T; L^2(\Omega))$, we have
\be
\int_{t_0}^{t_0+\tau}\int_{\Omega}{\widetilde{c_i}}^2 dxdt \le C^*\Gamma_4(1+\tau).
\la{widetildel2}
\ee
Using again (\ref{phinf}) we have
\be
\int_{t_0}^{t_0+\tau}\int_{\Omega}{{c_i}}^2 dxdt \le C^*\Gamma_5(1+\tau).
\la{ci2}
\ee
with constant $\Gamma^*_5$ depending like above only on $\Gamma$ and bounds on $\rho^*$. 

\vspace{.5cm}

\noindent{\bf{Step 4: Global bound on $c_i$ in $ L^{\infty}(L^2)$}.}\\
We use now (\ref{intfeq}) with $F(c) = \fr{c^2}{2}$. We have
\be
\fr{d}{dt}\int_{\Omega}c_i^2dx \le - 2D_i\int_{\Omega}|\na c_i|^2dx +
2D_i|z_i|\|c_i\|_{L^4(\Omega)}\|\na\Phi\|_{L^4(\Omega)}\|\na c_i\|_{L^2(\Omega)}.
\la{l2ineq}
\ee
We use the inequalities
\be
\|c_i\|_{L^4(\Omega)} \le C\left[\|\na c_i\|_{L^2(\Omega)}^{\fr{1}{2}} + \|c_i\|_{L^2(\Omega)}^{\fr{1}{2}}\right]\|c_i\|_{L^2(\Omega)}^{\fr{1}{2}}
\la{ladyb}
\ee
and we estimate
\be
\|\na\Phi\|_{L^4(\Omega)} \le \|\na(\Phi-\Phi^*)\|_{L^4(\Omega)} + \|\na\Phi^*\|_{L^4(\Omega)}.
\la{phiphistar}
\ee
For  $\|\na(\Phi-\Phi^*)\|_{L^4(\Omega)}$ we bound
\be
 \|\na(\Phi-\Phi^*)\|_{L^4(\Omega)} \le C \|\na(\Phi-\Phi^*)\|_{L^2(\Omega)}^{\fr{1}{2}}\|\rho-\rho^*\|_{L^2(\Omega)}^{\fr{1}{2}} \le C^*\Gamma^{\fr{1}{4}}\|\rho-\rho^*\|_{L^2(\Omega)}^{\fr{1}{2}}.
\la{phiphistarineq}
\ee
We used here that $\|\na(\Phi(t)-\Phi^*)\|_{L^2(\Omega)}^2$ is bounded in time because it is part of the energy. Putting these together we see that
\be
\ba
\fr{1}{2}\fr{d}{dt}\int_{\Omega}c_i^2dx \le -D_i\|\na c_i\|^2_{L^2(\Omega)} \\
+D_i\Gamma_6\|\na c_i\|_{L^2(\Omega)}\left[\|\na c_i\|_{L^2(\Omega)}^{\fr{1}{2}} + \|c_i\|_{L^2(\Omega)}^{\fr{1}{2}}\right]\|c_i\|_{L^2(\Omega)}^{\fr{1}{2}}\left[\sum_{j=1}^N\|c_j\|_{L^2(\Omega)}^{\fr{1}{2}} + \Gamma_7\right]
\ea
\la{l2sumi}
\ee
where the constants $\Gamma_6, \Gamma_7$ depend on the initial energy, $\epsilon$, all  $|z_j|$ and bounds on $\rho^*$, $\Phi^*$. From here we obtain
\be
\fr{d}{2dt} A^2 \le -\delta G^2 + \Gamma_6G(G^{\fr{1}{2}} + A^{\fr{1}{2}})A^{\fr{1}{2}}(A^{\fr{1}{2}} + \Gamma_7)
\la{ngineq}
\ee
for 
\be
A^2(t)= \sum_{j=1}^N\|c_j(t)\|^2_{L^2(\Omega)},\quad \quad G^2(t) = \sum_{j=1}^N\int_{\Omega}|\na c_j(x,t)|^2 dx,
\la{ag}
\ee
with slightly modified $\Gamma_6$ and $\Gamma_7$ and $\delta = \min{D_j}$. 
Using Young inequalities we finally obtain
\be
\fr{d}{dt} A^2 \le \Gamma_8(A^4 + A^2 ).
\la{aineq}
\ee
In view of (\ref{ci2}) we have 
\be
\int_{t_0}^{t_0+\tau} A^2 dt \le NC^*\Gamma_5(1+\tau),
\la{aint}
\ee
which, together with (\ref{aineq}) shows that $A$ remains bounded 
\be
\sup_{t_0\le t\le t_0+\tau}A^2(t) \le A(t_0)^2e^{\Gamma_9 (1+\tau)}
\la{l2pointwise}
\ee
where $\Gamma_9$ depends  on the initial energy, $\epsilon$, all  $|z_j|$ and bounds on $\rho^*$, $\Phi^*$. This is the first place where data appear in the right hand side of inequalities on their own and not through the initial energy $\Gamma$. Now we cover the interval $[0,T]$ with intervals of length $\fr{\tau}{2}$ where $\fr{\tau}{2}>0$ is a fixed positive time step. In view of (\ref{aint}) with $t_0 = 0$ and $\tau$ replaced by $\fr{\tau}{2}$, because of the Chebyshev inequality there exists $t_0\in [0,\fr{\tau}{2}]$ such that
\be
A(t_0)^2\le C^* \Gamma_5 \tau^{-1}.
\la{chebya}
\ee
Using this value we obtain from (\ref{l2pointwise})
\be
\sup_{\fr{\tau}{2} \le t\le \tau}A^2(t) \le  C^* \Gamma_5 e^{\Gamma_9 (1+\tau)}\tau^{-1}.
\la{l2pointwisen}
\ee
Now, because of (\ref{aint}) in the interval $[\fr{\tau}{2}, \tau]$ and the Chebyshev inequality, there is a new $t_0\in[ \fr{\tau}{2}, \tau]$ such that
(\ref{chebya}) holds, and thus, inductively
\be
\sup_{\fr{\tau}{2} \le t\le T}A^2(t) \le  C^* \Gamma_5 e^{\Gamma_9 (1+\tau)}\tau^{-1}.
\la{l2pointwisene}
\ee
This bound is independent of time, and depends only on initial energy and an arbitrary positive initial time $\fr{\tau}{2}>0$. We obtain also that $c_i\in L^{\infty}(0,T; L^2(\Omega))$, by adding the inequality (\ref{l2pointwise}) for the first time interval, starting at $t_0=0$ and obtain thus:
\be
\sup_{0 \le t\le T}A^2(t) \le  C^* (A(0)^2+\Gamma_5\tau^{-1}) e^{\Gamma_9 (1+\tau)} = \Gamma_{\tau}(1+A(0)^2).
\la{l2pointwisenew}
\ee
The right hand side does not depend of $T$. Returning to (\ref{ngineq})
we see that
\be
\sum_{i=1}^N\int_{t_0}^{t_0+\tau}\|\na c_i\|_{L^2(\Omega)}^2dt \le \Gamma_{\tau}(1+A(0)^2)
\la{gradcl2b}
\ee
with slightly different $\Gamma_\tau$. The mention of $\tau$ in the constant is only as a reminder of how the bound is achieved, but basically one thinks of $\tau =1$, i.e. a fixed auxilliary time step. 

\vspace{.5cm}

\noindent{\bf{Step 5: Global $L^{\infty}(L^p)$ bounds for $c_i$ and bounds for $\na\Phi$.}}\\
We improve the time integrability in (\ref{cilq1}) for $p>2$.  We write
\[
\int_{\Omega} c_i(x,t)^pdx = \int_{\Omega}c_i(x,t)^{2-\delta} c_i(x,t)^{p-2 +\delta}dx \le \left(\int_{\Omega}c_i(x,t)^2dx\right)^{1-\fr{\delta}{2}}\left(\int_{\Omega}c_i(x,t)^{\fr{2(p-2+\delta)}{\delta}}dx\right)^{\fr{\delta}{2}}
\]
and therefore, in view of (\ref{cilq1}) with $q = \fr{2(p-2+\delta)}{\delta}$ and (\ref{l2pointwisenew}), we have that
\be
\int_{t_0}^{t_0+\tau} \|c_i(t)\|_{L^p}^{\fr{p}{p-2  +\delta}}dt \le (\Gamma_{\tau}(1+A(0)^2))^\fr{2-\delta}{2(p-2+\delta)}\Gamma_3^*(1+\tau)
\la{intcipb}
\ee
holds for any $p>2$ and any $0<\delta<2$. 

By taking $2<p<4$ and $\delta$ small enough we have $\fr{p}{p-2+\delta}\ge 2$. Using the bound
\be
\|\na(\Phi-\Phi^*)\|_{L^{\infty}(\Omega)} \le C\|\rho-\rho^*\|_{L^p(\Omega)}
\la{gradphib}
\ee
we obtain that
\be
\int_{t_0}^{t_0+\tau} \|\na (\Phi-\Phi^*)\|_{L^{\infty}}^{2}dt \le \Gamma_{\tau}
\la{naphib}
\ee
holds with $\Gamma_{\tau}$ depending on initial energy, $\tau$ and $A(0)$. 
Using (\ref{intfineq}) with $F(c) = c^p$ and arbitrary $p\ge 2$ we obtain from (\ref{naphib})
\be
\sup_{0\le t\le T}\|c_i\|_{L^p(\Omega)} \le \Gamma_p
\la{linftylp}
\ee
with $\Gamma_p$ depending on initial energy, initial $\|c_i(0)\|_{L^p(\Omega)}$ but not on $T$. This is obtained in the same manner as the uniform bound 
(\ref{l2pointwisenew}): using controlled growth on overlapping short time intervals starting from values bounded using Chebyshev inequalities. Then, returning to the elliptic equation solved by $\Phi$ (\ref{poin}) we obtain uniform in time bounds for the norms of $\Phi$ in $W^{2,p}(\Omega)$. 
In particular,
\be
\|\Phi(\cdot, t)\|_{W^{1,\infty}(\Omega)} \le \Gamma_{\infty}^*
\la{phiginf}
\ee
holds for $t\ge 0$. 

\vspace{.5cm}

\noindent{\bf Step 6: Uniform bounds for $\c_i$.}\\
Now we turn to the equation satisfied by $\widetilde{c}_i$ 
\be
\pa_t{\widetilde{c_i}} = D_i\D {\widetilde{c_i}} - (u + D_iz_i\na\Phi)\na {\widetilde{c_i}} + z_i((\pa_t + u\cdot\na)\Phi){\widetilde{c_i}}.
\la{tildeciieq}
\ee
The boundary conditions are homogeneous Neumann:
\be
\pa_n\widetilde{c}_i(x,t)_{\left |\right. \;\pa\Omega} = 0.
\la{widetildebc}
\ee
Because of (\ref{gradcl2b}) and (\ref{phiginf}) we have that, for any $k=0, 1, \dots$ there exists
$t_k\in [k\fr{\tau}{2}, (k+1)\fr{\tau}{2}]$ such that
\be
\|\na\widetilde{c}_i(t_k)\|_{L^2(\Omega)}^2 \le \Gamma_{\tau}(1+ A(0))^2.
\la{tildeh1id}
\ee
If $c_i(0)\in H^1(\Omega)$ we can take $t_0=0$. We prove local uniform estimates 
\be
\sup_{t_k\le t\le t_k+\tau}\|\na\widetilde{c_i}(t)\|^2_{L^2(\Omega)} + \int_{t_k}^{t_k+ \tau}\|\D\widetilde{c}_i\|^2_{L^2(\Omega)}dt \le \Gamma_{\tau}.
\la{inth2tildecb}
\ee
These are obtained by multiplying (\ref{tildeciieq}) by $-\D\widetilde{c}_i$ and integrating. We obtain
\be
\ba
\fr{1}{2}\fr{d}{dt}\|\na\widetilde{c_i}(t)\|^2_{L^2(\Omega)} + D_i\|\D\widetilde{c}_i(t)\|^2_{L^2(\Omega)}\\
 \le C\left[\|u(t)\|_{L^4(\Omega)} +D_i|z_i|\|\na\Phi(t)\|_{L^4(\Omega)}\right]\|\na\widetilde{c_i}(t)\|^{\fr{1}{2}}_{L^2(\Omega)}\|\D\widetilde{c}_i(t)\|^{\fr{3}{2}}_{L^2(\Omega)}\\
  + |z_i|\|\pa_t\Phi(t) + u\cdot\na\Phi(t)\|_{L^4(\Omega)}\|\widetilde{c}_i(t)\|_{L^4(\Omega)}\|\D\widetilde{c}_i(t)\|_{L^2(\Omega)}.
\ea
\la{h1ineq}
\ee
Now we use a Gronwall inequality based on several facts. In view of (\ref{phiginf}) and the consequence
\be
\int_0^T\|u(t)\|_{L^4(\Omega)}^4dt\le \Gamma
\la{ul4}
\ee
of the energy inequality (\ref{gendiss}), the terms involving $u$ and $\na\Phi$ are easily bounded. The term involving $\pa_t\Phi$ is more interesting.
We use the Poisson equation and the equations (\ref{dtci}) to write
\be
\pa_t \Phi = \fr{1}{\epsilon}(-\D_D)^{-1}\left(\sum_{i=1}^N z_i\div\left(D_i\na c_i +(D_iz_i\na\Phi-u)c_i\right)\right).
\la{phiteq}
\ee
Because $(-\D_D)^{-1}\D$  is bounded in $L^4(\Omega)$ and $(-\D_D)^{-1}\div$ maps $L^2(\Omega)$ to $H_0^1(\Omega)\subset L^4(\Omega)$, we have
\be
\|\pa_t \Phi(t)\|_{L^4(\Omega)} \le C\sum_{i=1}^N\|c_i(t)\|_{L^4(\Omega)}(1+\|u(t)\|_{L^4(\Omega)} + \|\na\Phi(t)\|_{L^4(\Omega)}).
\la{phitl4}
\ee 
Because of (\ref{linftylp}) and (\ref{phiginf}), these  inequalities imply 
\be
\|\na\widetilde{c}_i(t)\|_{L^2(\Omega)}^2 + \int_{t_k}^t\|\D\widetilde{c}_i(s)\|^2_{L^2(\Omega)}ds \le \Gamma_{\tau}\|\na\widetilde{c}_i(t_k)\|_{L^2(\Omega)}^2 
\la{gronci}
\ee
for $t\in [t_k, t_k+\tau]$ and this implies (\ref{inth2tildecb}). Because $[(k+1)\fr{\tau}{2}, (k+2)\fr{\tau}{2}]\subset [t_k, t_k +\tau]$, from  (\ref{inth2tildecb}) we deduce by induction
\be
\sup_{0\le t\le T}\|\na\widetilde{c_i}(t)\|^2_{L^2(\Omega)}\le \Gamma_{\tau}
\la{tildecih1bound}
\ee
and
\be
\int_0^T\|\D\widetilde{c}_i(t)\|^2_{L^2(\Omega)}dt \le \Gamma_{\tau}T.
\la{tildeh2b}
\ee
Returning to the local estimates, we find new $t_k\in [k\fr{\tau}{2}, (k+1)\fr{\tau}{2}]$ such that
\be
\|\D\widetilde{c}_i(t_k)\|_{L^2(\Omega)}^2 \le\Gamma_{\tau}
\la{h2points}
\ee
for $k\ge 0$. We use now a local energy estimate for the Navier-Stokes equation:
\be
\sup_{t_k\le t\le t_k+\tau}\|\na u(t)\|^2_{L^2(\Omega)} + \nu\int_{t_k}^{t_k+ \tau}\|\D u(t)\|_{L^2(\Omega)}^2dt \le \Gamma_{\tau}
\la{nsestr}
\ee
which is based on the fact that the forcing in (\ref{nse}) is bounded in $L^2(\Omega)$ and on standard estimates for the nonlinearity and the Stokes operator.
Using the embedding $H^2(\Omega)\subset L^{\infty}(\Omega)$ we have thus
\be
\int_{t_k}^{t_k + \tau}\|u(t)\|_{L^{\infty}(\Omega)}^2 dt \le \Gamma_{\tau}
\la{uinfty}
\ee
and, from (\ref{inth2tildecb}),
\be
\int_{t_k}^{t_k + \tau}\|\widetilde{c_i}(t)\|_{L^{\infty}(\Omega)}^2 dt \le \Gamma_{\tau}.
\la{tildeciinfty}
\ee
Now we take a large $p$ and estimate from (\ref{tildeciieq})
\be
\ba
\fr{1}{p}\fr{d}{dt}\|\widetilde c_i(t)\|_{L^p(\Omega)}^p + (p-1)D_i\int_{\Omega}|\na\widetilde{c}_i(x,t)|^2\widetilde{c}_i(x,t)^{p-2}dx\\
\le D_i|z_i|\|\na\Phi(t)\|_{L^{\infty}(\Omega)}\int_{\Omega}|\na\widetilde{c}_i(x,t)|\widetilde{c}_i(x,t)^{p-1}dx + |z_i|\|\pa_t\Phi + u\na\Phi\|_{L^{\infty}(\Omega)}\|\widetilde c_i(t)\|_{L^p(\Omega)}^p.
\ea
\la{lpwidetildecib}
\ee
Consequently, 
\be
\sup_{t_k\le t\le t_k + \tau}\|\widetilde c_i(t)\|_{L^p(\Omega)} \le\|\widetilde c_i(t_k)\|_{L^p(\Omega)}e^{\int_{t_k}^{t_k+\tau}\left[D_i\fr{1}{2(p-1)}|z_i|^2\|\na\Phi(t)\|_{L^{\infty}(\Omega)}^2 + |z_i|\|\pa_t\Phi(t) + u(t)\na\Phi(t)\|_{L^{\infty}(\Omega)}\right]dt}.
\la{lpwidetildecibound}
\ee
Passing  $p\to\infty$ we have
\be
\sup_{t_k\le t\le t_k + \tau}\|\widetilde c_i(t)\|_{L^\infty(\Omega)}
\le \|\widetilde c_i(t_k)\|_{L^{\infty}(\Omega)}e^{\int_{t_k}^{t_k+\tau}|z_i|\|\pa_t\Phi(t) + u(t)\na\Phi(t)\|_{L^{\infty}(\Omega)}dt}.
\la{tildelinftyloc}
\ee  
Using  (\ref{phiteq}) we have now enough information to bound,
\be
\|\pa_t\Phi(t)\|_{L^{\infty}(\Omega)}\le C\sum_{i=1}^N[\|\na c_i\|_{L^p(\Omega)}+ \|D_iz_i\na\Phi-u\|_{L^{\infty}(\Omega)}\|\na c_i\|_{L^2(\Omega)} +\|\rho c_i\|_{L^2(\Omega)}]
\la{patphiineq}
\ee
with $p>2$, where we also used the fact that $(-\D_D)^{-1}\div$ maps $L^p(\Omega)$ to $L^{\infty}(\Omega)$ and $(-\D_D)^{-1}$ maps $L^2(\Omega)$ to $L^{\infty}(\Omega)$. 
Because $\na c_i = e^{-z_i\Phi}(\na\widetilde{c_i} -z_i\widetilde{c_i}\na\Phi)$, the bound (\ref{phiginf}), the embedding $H^2(\Omega)\subset W^{1,p}(\Omega)$ and (\ref{inth2tildecb}) we have
\be
\int_{t_k}^{t_k+\tau}\|\pa_t\Phi(t)\|_{L^{\infty}(\Omega)}dt \le \Gamma_{\tau}
\la{intphit}
\ee
and consequently, by induction, we obtain the uniform bound
\be
\sup_{0\le t\le T}\|\widetilde{c_i}(t) \|_{L^{\infty}(\Omega)}\le \Gamma_{\tau}.
\la{tildecinfty}
\ee
This then implies
\be
\sup_{0\le t\le T}\|{c_i}(t) \|_{L^{\infty}(\Omega)}\le \Gamma_{\tau}.
\la{cinfty}
\ee
In view of (\ref{natildeci}) we have using (\ref{cinfty})
\be
\int_0^T\|\na\widetilde{c}_i(t)\|_{L^2(\Omega)}^2dt \le \Gamma
\la{globdiss}
\ee
with $\Gamma$ time independent. The Nernst-Planck equations (\ref{cpmeq}) imply
\be
\pa_t\rho = \sum_{i=1}^Nz_iD_i\div\left(e^{-z_i\Phi}\na\c_i\right) -u\cdot\na\rho
\la{rhoteq}
\ee
and, together with the Poisson equation (\ref{poiphi}), the bound (\ref{phinf}) and the embedding $H^1(\Omega)\subset L^p(\Omega)$  have then the consequence that
\be
\int_0^T\|\pa_t\Phi(t)\|_{L^p(\Omega)}^2dt \le \Gamma
\la{phitlp}
\ee
holds for any $p\in [1,\infty)$. The $L^2$ boundedness of Dirichlet Riesz transforms imply that 
\be
\int_0^T\|\na\pa_t\Phi(t)\|_{L^2(\Omega)}^2dt \le \Gamma
\la{naphitl2}
\ee
also holds. Turning to the equation (\ref{tildeciieq}) we obtain
\be
\la{deltildec}
\int_0^T\left [\|\pa_t\c_i(t)\|_{L^2(\Omega)}^2 + \|\D\c_i(t)\|_{L^2(\Omega)}^2\right]dt \le \Gamma.
\ee
Indeed, 
\be
\int_0^T \left[\|u(t)\|_{L^4(\Omega)}^4 +\|\na\Phi(t)\|_{L^{\infty}(\Omega)}^2\right]\|\na\c_i(t)\|_{L^2(\Omega)}^2dt\le \Gamma
\la{coef1bound}
\ee
and 
\be
\int_0^T \|(\pa_t\Phi(t) + u(t)\cdot\na\Phi(t))\c_i(t)\|_{L^2(\Omega)}^2dt\le \Gamma
\la{coef2bound}
\ee
because of  (\ref{phiginf}), (\ref{ul4}), (\ref{tildecih1bound}), (\ref{tildecinfty}), (\ref{globdiss}) and (\ref{phitlp}).

\beg{thm}\la{globlock}
Let $\Omega\subset \Rr^2$  be a bounded domain with smooth boundary. Let $z_i\in\Rr$, $1\le i\le N$, let $\epsilon>0$, $D_i>0$, $i=1,\dots N$. Let $c_i(0)$ be nonnegative functions $1\le i\le N$, with $c_i(0)\in L^p(\Omega)\cap W^{2,q}(\Omega)$, $p=2q>4$, let $W\in W^{\fr{3}{2}, p}(\pa\Omega)$ be a  function defined on $\pa\Omega$ and let $u_0\in W_0^{1,p}(\Omega)$ be divergence-free. Then there exists a unique global strong solution with initial data $c_i(0)$,
$u_0$, defined on $[0,\infty)$ of the Nernst-Planck-Navier-Stokes system
\be
\left \{
\ba
\pa_t u + u\cdot\na u + \na p = \nu\D u -  (k_B\T) \rho\na\Phi,\\
\div u = 0,\\  
\rho = \sum_{i=1}^N z_i c_i,\\
-\epsilon\D \Phi = \rho,\\
\pa_t c_i + u\cdot\na c_i = D_i\left(\D c_i + z_i{\mbox{div}}(c_i\na\Phi)\right),
\ea
\right .
\la{pnpns}
\ee
in $\Omega\times [0,\infty)$, with boundary conditions
\be
\left\{
\ba
u_{\left |\right. \pa \Omega} = 0,\\
\Phi_{\left |\right. \pa \Omega} = W,\\
(\na c_i + z_ic_i\na\Phi)_{\left |\right. \pa \Omega}\;\cdot n = 0
\ea
\right.
\la{pnpnsbc}
\ee
where $n$ is the external normal at $\pa\Omega$. There exist constants $\Gamma_p$ depending on the parameters $\epsilon, D_i, z_i$, the domain $\Omega$, the initial energy ${\mathcal E(0)}$, and the norms 
\be
\|c_i(0)\|_{L^p(\Omega)}, \quad  \|W\|_{H^s(\pa\Omega)},\quad  \|u_0\|_{L^2(\Omega)},
\la{norms}
\ee
for $p\ge 2$, $s\ge \fr{3}{2}$, such that 
\be
\max_{1\le i\le N}\sup_{0\le t<\infty }\|c_i(t)\|_{L^{p}(\Omega)} \le \Gamma_{p},
\la{linftylpa}
\ee
The bounds 
\be
\sup_{0\le t<\infty }\|\Phi(t)\|_{W^{2,p}(\Omega)}\le \Gamma_p.
\la{phi2wbound}
\ee
hold for $p\ge 2$ and in particular, 
\be
\sup_{0\le t<\infty}\|\Phi(t)\|_{W^{1,\infty}(\Omega)} \le \Gamma_{\infty}^*
\la{phinfa}
\ee
holds.  In addition 
\be
\max_{1\le i\le N}\sup_{0\le t<\infty}\|c_i(t)\|_{L^{\infty}(\Omega)} \le \Gamma_{\infty},
\la{linftyinfty}
\ee
and
\be
\max_{1\le i\le N}\sup_{0\le t<\infty }\|\na c_i(t)\|_{L^2(\Omega)}\le \Gamma.
\la{h1cbound}
\ee
The functions
\[
\widetilde{c_i}(x,t) = c_i(x,t)e^{z_i\Phi(x,t)}
\]
defined in (\ref{tildeci}) obey
\be
\max_{1\le i\le N}\sup_{0\le t<\infty}\|\widetilde{c}_i(t)\|_{L^{\infty}(\Omega)} \le \Gamma_{\infty},
\la{linftytildeinfty}
\ee
\be
\max_{1\le i\le N}\sup_{0\le t<\infty}\|\na \widetilde{c}_i(t)\|^2_{L^2(\Omega)}\le \Gamma,
\la{h1ctildebound}
\ee
and
\be
\int_0^{\infty}\int_{\Omega} \left |\na{\widetilde{c_i}}(x,t)\right|^2 dxdt \le \Gamma_2.
\la{nawdtld}
\ee
Moreover,
\be
\int_0^{\infty} \left[ \|\pa_t\widetilde{c}_i(t)\|_{L^2(\Omega)}^2 + \|\D\widetilde{c}_i(t)\|_{L^2(\Omega)}^2\right]dt \le \Gamma
\la{tildetimeint}
\ee
and
\be
\int_0^{\infty}\left[ \|\pa_t\Phi(t)\|_{L^p(\Omega)}^2 + \|\na\pa_t\Phi(t)\|^2_{L^2(\Omega)}\right]dt\le \Gamma
\la{patphinaphi}
\ee
hold. The Navier-Stokes solution satisfies 
\be
\sup_{0\le t\le T}\|u(t)\|_{H^1(\Omega)}^2 + \nu \int_0^T\|u(t)\|^2_{H^2(\Omega)}dt \le C^*T
\la{strongns}
\ee
for any $T>0$, with $C^*$ depending on $\Gamma_p$ above and $\|u_0\|_{H^1(\Omega)}$
and further
\be
\sup_{0\le t\le T}[\|c_i(t)\|_{W^{2,q}(\Omega)} + \|\pa_tc_i(t)\|_{L^q(\Omega)} ]\le C_q(T)
\la{maxlqT}
\ee
and
\be
\int_0^{T}[\|\pa_t u(t)\|^2_{L^p(\Omega)} + \|u(t)\|_{W^{2,p}(\Omega)}^2]dt
\le U_p(T)
\la{maxregat}
\ee
hold. The constants $C_q(T)$ and $U_p(T)$ depend on the initial data and $T$.
\end{thm} 
The proof of the theorem follows immediately from the a priori bounds established above and a uniform local existence and uniqueness theorem.
\beg{thm}\la{loc} Let $\Omega\subset \Rr^d$, $d=2,3$ be a bounded domain with smooth boundary. Let $z_i\in\Rr$, $1\le i\le N$ and let $\epsilon>0$, $D_i>0$, $i=1,\dots N$. Let $c_i(0)$ be nonnegative 
functions $1\le i\le N$, let $W$ be a smooth function defined on $\pa\Omega$  and let $u_0\in H^1(\Omega)^d$ be divergence-free. Let $p=2q>2d$. There exists $T_0$ depending only on the parameters of the problem $\epsilon, D_i, z_i, \nu$, the domain $\Omega$, the initial energy $\mathcal E(0)$ and on the norms
\be
\|c_i(0)\|_{L^p(\Omega)}, \quad  \|W\|_{W^{\fr{3}{2},p}(\pa\Omega)}, \quad \|u_0\|_{L^{2}(\Omega)},
\la{normssa}
\ee
such that a unique strong solution of (\ref{pnpns}) with initial data 
$c_i(0)\in L^p(\Omega)\cap W^{2,q}(\Omega)$, $u_0\in W^{1,p}(\Omega)$ exists and satisfies 
\be
\sup_{0\le t\le T_0}\|c_i(t)\|_{L^p(\Omega)} \le 3 \|c_i(0)\|_{L^p(\Omega)}
\la{cilptriple}
\ee
and 
\be
\sup_{0\le t\le T_0}[\|c_i(t)\|_{W^{2,q}(\Omega)} + \|\pa_tc_i(t)\|_{L^q(\Omega)}] \le C_q
\la{maxlq}
\ee
and
\be
\int_0^{T_0}[\|\pa_t u(t)\|^2_{L^p(\Omega)} + \|u(t)\|_{W^{2,p}(\Omega)}^2]dt 
\le C_p
\la{maxrega}
\ee
with constants $C_p, C_q$, depending on
\be
\|c_i(0)\|_{L^p(\Omega)}, \|c_i(0)\|_{W^{2,q}(\Omega)}, \|u(0)\|_{W^{1,p}(\Omega)}.
\la{initialnorms}
\ee
\end{thm}
\beg{rem} Note that the time of existence depends only on the initial energy and the norms of $c_i(0)$, $u_0$, $W$ listed in (\ref{normssa}), but not on the higher norms which are subsequently controlled. There is no special meaning to the time $t=0$: the result holds from any $t_0$ for a short time, determined as above. We also remark that although no attempt was made to find the most generous initial data regularity conditions, nevertheless no compatibility conditions for the initial data are  required.
\end{rem}
The proof is presented in Appendix B. Using the global existence theorem we obtain
\beg{thm}\la{decay} Under the conditions of Theorem \ref{globlock} we obtain
\be
\lim_{t\to\infty}\int_{\Omega}\left |\na\widetilde{c}_i(x,t)\right|^2dx = 0.
\la{nactwidletozero}
\ee
Consequently we also have
\be
\lim_{t\to\infty}{\mathcal {D}}(t) = 0
\la{dtozero}
\ee
where
\be
{\mathcal{D}}(t) = \sum_{i=1}^N D_i\int_{\Omega}c_i(x,t)\left |\na\left(\log\left(c_i(x,t)e^{z_i\Phi(x,t)}\right)\right)\right|^2 dx.
\la{mathcalditself}
\ee
\end{thm} 
\noindent{\bf{Proof.}}
Let 
\be
N(t) = \|\na \c_i(t)\|^2_{L^2(\Omega)}.
\la{nat}
\ee
The proof of Theorem \ref{decay} is done by contradiction. Assume by contradiction that there exists a sequence of times $t_n\to\infty$ where
\be
N(t_n) = \|\na\widetilde{c}_i(t_n)\|_{L^2(\Omega)}^2\ge \delta>0.
\la{absurdass}
\ee
The time derivative of $N(t)$ is
\be
N'(t) = 2\int_{\Omega}\pa_t{\widetilde c_i}(x,t)(-\D\widetilde{c}_i(x,t))dx.
\la{nprime}
\ee
In view of (\ref{tildetimeint}) we have that
\be
\int_0^{\infty}|N'(t)|dt \le \Gamma<\infty.
\la{nprimel1}
\ee
Thus the limit 
\be
N(\infty)= \lim_{t\to\infty} N(t) = N(0) + \int_0^{\infty}N'(t)dt
\la{limN}
\ee
exists, and by the contradiction assumption $N(\infty)\ge \delta>0$. Therefore there exists $T>0$ such that $N(t)\ge \fr{\delta}{2}$ for all $t\ge T$. This is absurd, because
\be
\int_0^{\infty}N(t)dt \le \Gamma_2 <\infty
\la{intN}
\ee
by (\ref{nawdtld}).

We prove now convergence of solutions for infinite time. 
\beg{thm}\la{asy}
Let the conditions of Theorem \ref{globlock} be satisfied. Then the solution converges to a Boltzmann state, and the velocity converges to zero. The Boltzmann state is uniquely determined by the initial concentrations
\be
I_i^0 = \int_{\Omega}c_i(0)dx,
\la{iizero}
\ee
and has the form
\be
c_i^* = Z_i^{-1}e^{-z_i\Phi^*}
\la{bof}
\ee
with 
\be
Z_i = (I_i^0)^{-1}\int_{\Omega}e^{-z_i\Phi^*}dx
\la{ziphinfty}
\ee
and with $\Phi^*$ solving  
\be
-\epsilon\D\Phi^*= \sum_{i=1}^N z_i I_i^0 \fr{e^{-z_i\Phi^*}}{\int_{\Omega} e^{-z_i\Phi^*}dx}
\la{phinftyeq}
\ee
with boundary conditions (\ref{phibc}). 
\end{thm}
\noindent{\bf Proof.} 
Because of the boundary condition $\pa_n\widetilde{c}_i =0$, we have that
\be
\|\widetilde{c}_i (\cdot,t)-m_i(t)\|_{L^2(\Omega)} \le C\|\na\widetilde{c}_i\|_{L^2(\Omega)}
\la{poinca}
\ee
where
\be
m_i(t) = |\Omega|^{-1}\int_{\Omega}\widetilde{c}_i(x,t)dx.
\la{mi}
\ee
Thus, from (\ref{nawdtld}) we have
\be
\int_0^{\infty}\|\widetilde{c}_i (\cdot,t)-m_i(t)\|_{L^2(\Omega)}^2dt \le \Gamma
\la{cimil2}
\ee
and, in view of the convergence
\be
\lim_{t\to\infty}\|\na\widetilde{c}_i(t)\|_{L^2(\Omega)} = 0
\la{widetildel2tozero}
\ee
and the boundary condition $\pa_n\widetilde{c}_i =0$, we have that
\be
\lim_{t\to\infty}\|\widetilde{c}_i (\cdot,t)-m_i(t)\|_{L^2(\Omega)} = 0.
\la{limtildem}
\ee
Let $s_n\to\infty$ be any sequence of times. By extracting a subsequence denoted $t_n$, in view of the previous results, we may assume without loss of generality that there exist numbers
$M_i\ge 0$ and a function $\Phi_{\infty}$ such that
\be
\lim_{n\to\infty}m_i(t_n)  = M_i,
\la{mitnm}
\ee
\be
\lim_{n\to\infty}\|\widetilde{c}_i(t_n) -M_i\|_{H^1(\Omega)} = 0
\la{convtilde}
\ee
holds in $H^1(\Omega)$ and 
\be
\lim_{n\to\infty} \|\Phi(t_n)-\Phi_{\infty}\|_{W^{1,\infty}(\Omega)} = 0
\la{phiconv}
\ee
holds in $W^{1,\infty}(\Omega)$ by compactness of the embedding $W^{2,p}(\Omega)\subset \subset W^{1,\infty}(\Omega)$ for $p>2$. Then it follows from the above that
\be
\lim_{n\to\infty}\|c_i(t_n) - M_ie^{-z_i\Phi_{\infty}}\|_{H^1(\Omega)} = 0.
\la{ciconv}
\ee
In addition,
\be
\int_{\Omega}c_i(x,t_n)dx = I_i^0 
\la{izero}
\ee
follows from the zero flux boundary conditions, and thus we identify the constants $M_i$ as
\be
M_i = I_i^0\left(\int_{\Omega} e^{-z_i\Phi_{\infty}}dx\right)^{-1}.
\la{mid}
\ee
Passing to the limit in the equation (\ref{poin}) we have therefore that
$\Phi_{\infty}$ solves (\ref{phinftyeq}).
We remark that this equation does not depend on the sequence $s_n$. The proof of Theorem \ref{asy} is completed by the uniqueness of solutions of (\ref{phinftyeq}), proved below.

\beg{thm}\la{uniqphinfty} Any two $W^{1,\infty}(\Omega)$ solutions of (\ref{phinftyeq}) with the same Dirchlet boundary conditions (\ref{phibc}) must coincide.
\end{thm}

\noindent{\bf{Proof}}. Indeed, let $\Phi_{\infty}^{(i)}$, $i=1,2$ be the two solutions and let $\psi$ be their difference,
\be
\psi = \Phi_{\infty}^{(2)}-\Phi_{\infty}^{(1)}.
\la{psiphi}
\ee
Then $\psi$ satisfies
\be
-\epsilon\D\psi(x) + \sum_{i=1}^Nz_i^2I_i^0\int_0^1 p^i_{\lambda}(x)\left(\psi(x)- (\psi, p^i_{\lambda})_{L^2(\Omega)}\right)d\lambda = 0
\la{psieq}
\ee
with homogeneous boundary conditions. Here
\be
p^i_{\lambda} = \fr{e^{-z_i\Phi_{\lambda}}}{\int_{\Omega}e^{-z_i\Phi_{\lambda}}dx}
\la{pil}
\ee
and
\be
\Phi_{\lambda} = \Phi_{\infty}^{(1)} + \lambda \psi.
\ee
Taking the scalar product of (\ref{psieq}) with $\psi$ we obtain
\be
0 = \epsilon \|\na\psi\|^2_{L^2(\Omega)} + \sum_{i=1}^N z_i^2I_i^0\int_0^1d\lambda \int_{\Omega}p^i_{\lambda}(x)\left(\psi(x)- (\psi, p^i_{\lambda})_{L^2(\Omega)}\right)^2dx
\la{cov}
\ee
and therefore $\psi=0$. This uses the fact that $p^i_{\lambda}$ are probabilities. The existence of solutions of (\ref{phinftyeq}) is classical \cite{friedman}. We briefly discuss the Poisson-Boltzmann equations encountered in the present work, including the nice structure of their linearizations in Appendix A.

\section{Unconditional global stability for uniform selective boundary conditions}
In this section we consider uniform selective boundary conditions (\ref{gammaibc}). We remark that we only use the uniform aspect, i.e. the constancy of $\gamma_i$ and of $w_i = W(x)_{\left|\right. S_i}$, for the decay in Theorem \ref{diss}.

The positivity of $c_i(x,t)$  follows in exactly the same way as in the case of blocking boundary conditions: the equation (\ref{intfeq}) holds because, for $i\le M$ and $x\in S_i$ we have that $F'(\gamma_i) = 0$ and for $x\in \pa\Omega\setminus S_i$ the normal flux vanishes, and thus integration by parts is allowed. The steps 1, 2, and 3 of the proof for blocking boundary conditions are still valid: they do not use boundary conditions for $c_i$. In particular (\ref{phinf}),
(\ref{cilq1}), and (\ref{ci2})  still hold.

\vspace{.5cm}

\noindent{\bf{ Step 4: Global bound on $c_i$ in $L^{\infty}(L^2)$.}}\\
We introduce smooth time independent functions $g_i$ for $i=1, \dots, M$ such that
\be
{g_i}_{\left |\right. \; S_i} = \gamma_i.
\la{gigammai}
\ee
The evolution equations (\ref{cpmeq}) can be written as
\be
\pa_t c_i = D_i\div\left(\na c_i + z_ic_i\na\Phi \right) - u\cdot\na c_i.
\la{cifluxeq}
\ee
Multiplying by $c_i-g_i$ and using the boundary conditions (\ref{gammaibc}) which imply that
\be
(c_i-g_i)(\pa_n c_i + z_ic_i\pa_n\Phi)_{\left |\right. \pa\Omega} = 0,\quad {\mbox{for}}\;  i=1,\dots N,
\la{bcgi}
\ee
we obtain after integration by parts 
\be
\ba
\fr{1}{2}\fr{d}{dt}\int_{\Omega}(c_i^2(x,t)-2 g_i(x)c_i(x,t))dx + D_i\int_{\Omega}|\na c_i(x,t)|^2dx \\
\le D_i|z_i|\|c_i\|_{L^4(\Omega)}\|\na\Phi\|_{L^4(\Omega)}\left[\|\na c_i\|_{L^2(\Omega)} + \|\na g_i\|_{L^2(\Omega)}\right] + D_i\|\na c_i\|_{L^2(\Omega)}  \|\na g_i\|_{L^2(\Omega)}.
\ea
\la{cl2ineq}
\ee
Because $g_i, \na g_i$ are bounded, and the inequality (\ref{phiphistarineq}) is still valid,  the quantity
\be
y(t) = \sum_{i=1}^N \int_{\Omega}(c_i^2(x,t)-2 g_i(x)c_i(x,t))dx
\la{yt}
\ee
obeys 
\be
\fr{dy}{dt} \le C(y^2+1)
\la{ytineq}
\ee
and the local uniform bound $\int_{t_0}^{t_0+\tau} y(t)dt \le \Gamma (1+\tau)$.
Using a similar argument as in Step 4 of the proof of Theorem \ref{globlock} we deduce  the inequalities (\ref{l2pointwisenew}) and
(\ref{gradcl2b}).

\vspace{.5cm}

\noindent{\bf Step 5: Global $L^{\infty}(L^p)$ bounds for $c_i$ and bounds for $\na\Phi$.}\\
The inequality (\ref{naphib}) is obtained without use of boundary conditions for $c_i$ from (\ref{l2pointwisenew}), in the same manner as for the blocking boundary conditions case, and so it is thus still valid. We take the equations (\ref{cifluxeq}), multiply by $F'(c_i)- F'(g_i)$ where $F(c) = c^p$ and integrate by parts.  The boundary terms vanish, and thus, after integrating by parts we obtain
\be
\fr{d}{dt}\int (F(c_i(t)) - c_i(t)F'(g_i))dx  = -D_i\int_{\Omega} \left(\na c_i + z_i c_i\na\Phi\right)\na (F'(c_i)-F'(g_i))dx - \int_{\Omega}c_i u\cdot\na ( F'(g_i))dx.
\la{lpintev}
\ee
Using (\ref{naphib}) and (\ref{gradcl2b}) we obtain like in the case of blocking boundary conditions (\ref{linftylp})
\be
\sup_{0\le t}\|c_i(t)\|_{L^p(\Omega)}\le \Gamma_p
\la{cip}
\ee
and consequently
\be
\sup_{0\le t}\|\Phi (t)\|_{W^{2,p}(\Omega)}\le \Gamma_p
\la{phiw2p}
\ee
with $\Gamma_p$ time independent.

\vspace{.5cm}
\noindent{\bf{Step 6: Uniform bounds for $\c_i$.}}\\
These are obtained in the exact same manner as in the case of blocking boundary conditions. The auxilliary functions $\c_i$ obey time independent Dirichlet boundary conditions on $S_i$, for $i\le M$ and homogeneous Neumann conditions on the rest of the boundary and for $i\ge M+1$. Therefore 
\be
(\pa_t\c_i\pa_n \c_i)_{\left|\right.\; \pa\Omega} = 0, \quad i=1,\dots, N,
\la{seltildebc}
\ee
and thus there is no contribution from the boundary when we multiply the equation obeyed by $\c_i$ (\ref{tildeciieq}) by $-\D\c_i$ and integrate.
The rest of the arguments are repeated almost verbatim. We have thus
\beg{thm}\la{uniblock} Under the assumptions of Theorem \ref{globlock} a unique global strong solution of the Nernst-Planck-Navier-Stokes system (\ref{pnpns})  with uniform selective boundary conditions 
\be
\left\{
\ba
u_{\left |\right. \pa\Omega}=0,\\
\Phi_{\left |\right. \pa \Omega} = W(x),\\
W(x)_{\left |\right. S_i}= w_i, \quad i=1, \dots, M,\\
c_i{\left |\right. S_i}= \gamma_i, \quad i=1, \dots, M,\\
(\pa_nc_i + z_ic_i\pa_n\Phi)_{\left |\right. \pa \Omega\setminus S_i} = 0,\quad i= 1,\dots, M,\\
(\pa_nc_i + z_ic_i\pa_n\Phi)_{\left |\right. \pa \Omega} = 0,\quad i= M+1,\dots, N
\ea
\right.
\la{uselbc}
\ee
exists. The solution obeys the inequalities (\ref{linftylpa})-(\ref{maxregat}). As time tends to infinity, the velocity tends to zero and the solutions $c_i$ converge to the Boltzmann state 
\be
c_i^* = Z_i^{-1}e^{-z_i\Phi^*}
\la{bs}
\ee
with
$Z_i$ given by (\ref{zichoice}) for $i=1, \dots M$ and (\ref{ziphinfty}) for $i= M+1,\dots, N$ where $I_i^0$ are given by (\ref{iizero}), and with $\Phi^*$ solving the Poisson-Boltzmann equation
\be
-\epsilon \D\Phi^* = \rho^* = \sum_{i=1}^Nz_ic_i^*
\la{poiusel}
\ee
with boundary conditions (\ref{phibc}).
\end{thm}
\section{Global existence for general selective boundary conditions}
The case of general selective boundary conditions is different because the decay in Theorem \ref{diss} is no longer generally true. We can however use the
dissipative structure to obtain time dependent bounds, which allow for growth of norms but no finite time singularities. The approach is similar to the one for blocking and uniform selective boundary conditions once the replacement of the first step is obtained. This is done as follows. We start from the fundamental structure (\ref{dter}) for the energy density (\ref{densen}) relative to a Boltzmann state with $Z_i>0$ chosen below.
In view of (\ref{deleicin}) we observe that the general selective boundary conditions imply that the densities 
\be
\left(\fr{\delta \mathcal E}{\delta c_i}\right)_{\left |\right. S_i} = \log \gamma_i + z_i W(x) + \log Z_i
\la{delebc}
\ee
are known on the boundary for $i=1, \dots, M$. We consider a smooth, time independent function $\widetilde{W}(x)$ of $x\in\Omega$ such that
\be
\widetilde{W}(x)_{\left |\right. S_i} = W(x)
\la{widetildeW}
\ee
and choose
\be
Z_i = (\gamma_i)^{-1}
\la{zichoicen}
\ee
for $i=1, \dots, M.$ The rest of $Z_i$ may be arbitrary positive numbers.
We then write (\ref{dter}) as
\be
D_t E = \sum_{i=1}^N \left(\fr{\delta\mathcal E}{\delta c_i}-z_i\widetilde{W} +z_i\widetilde{W}\right)D_tc_i  - F\cdot u + R,
\la{dterhi}
\ee
add $\fr{|u|^2}{2k_BT_K}$ and integrate. Moving the term $z_i\widetilde{W}\pa_t c_i$ to the left hand side, integrating by parts using (\ref{cieq}) and the fact that on the selected portions $S_i$ of the boundary we have that $\fr{\delta\mathcal E}{\delta c_i}-z_i\widetilde{W}=0$, while on all the rest the normal derivative $\pa_n\fr{\delta\mathcal E}{\delta c_i} =0$, we obtain
\be
\ba
\fr{d}{dt}\left[{\mathcal{E}} + \fr{\|u(t)\|^2_{L^2(\Omega)}}{2k_BT_K} -
\sum_{i=1}^N\int_{\Omega}z_i\widetilde{W}(x)c_i(x,t)dx\right] + \mathcal D + \fr{\nu}{k_BT_K}\|\na u(t)\|^2_{L^2(\Omega)}\\
= -\sum_{i=1}^N\int_{\Omega}z_ic_i(x,t)u(x,t)\cdot\na \widetilde{W}(x)dx +\sum_{i=1}^ND_i\int_{\Omega}z_ic_i\na\left(\fr{\delta\mathcal E}{\delta c_i}\right)\na {\widetilde{W}}dx.
\ea
\la{dishi}
\ee
Thus
\be
\ba
\fr{d}{dt}\left[{\mathcal{E}} + \fr{\|u(t)\|^2_{L^2(\Omega)}}{2k_BT_K} -
\int_{\Omega}\rho(x,t)\widetilde{W}(x)dx\right] + \mathcal D + \fr{\nu}{k_BT_K}\|\na u(t)\|^2_{L^2(\Omega)}\\
= -\int_{\Omega}\rho(x,t)u(x,t)\cdot\na \widetilde{W}(x)dx +\sum_{i=1}^ND_iz_i\int_{\Omega}c_i\na\left(\fr{\delta\mathcal E}{\delta c_i}\right)\na {\widetilde{W}}dx.
\ea
\la{disw}
\ee
Let
\be
\mathcal F = {\mathcal{E}} + \fr{\|u(t)\|^2_{L^2(\Omega)}}{2k_BT_K} -
\int_{\Omega}\rho(x,t)\widetilde{W}(x)dx,
\la{calf}
\ee
\be
\mathcal G = {\mathcal{E}} + \fr{\|u(t)\|^2_{L^2(\Omega)}}{2k_BT_K},
\la{calg}
\ee
and
\be
{\mathcal D_1} = \mathcal D + \fr{\nu}{k_BT_K}\|\na u(t)\|^2_{L^2(\Omega)}.  
\la{done}
\ee
We remark that
\be 
|\mathcal F -\mathcal G| \le C(1 + \mathcal E)
\la{fgineq}
\ee
because $\widetilde{W}$ is bounded.
The first term in the right hand side of (\ref{disw}) can be estimated as follows,
\be
\ba
\left |\int_{\Omega}\rho(x,t)u(x,t)\cdot\na \widetilde{W}(x)dx\right| =
\left|\int_{\Omega}(\rho(x,t)-\rho^*(x,t) + \rho^*(x,t))u(x,t)\cdot\na \widetilde{W}(x)dx\right|\\
\le \left|\epsilon\int_{\Omega}\D(\Phi-\Phi^*)u\cdot\na\widetilde{W}dx\right| +
\int_{\Omega}|\rho^*||u(x,t)||\na\widetilde{W}|dx
\le \epsilon\int_{\Omega}|\na(\Phi-\Phi^*)|\na(u\cdot\na{\widetilde{W}})|dx + C\sqrt{\mathcal G}\\
\le C{\mathcal G} + C(1+\sqrt{\mathcal D_1}){\sqrt{\mathcal G}}.
\ea
\la{rhsdis1}
\ee
The second term is estimated using the dissipation $\mathcal D$ and the boundedness of $\na\widetilde{W}$,
\be
\ba
\left |\sum_{i=1}^ND_iz_i\int_{\Omega}c_i\na\left(\fr{\delta\mathcal E}{\delta c_i}\right)\na {\widetilde{W}}dx\right| \le
C\sqrt{\mathcal D}\sqrt{\sum_{i=1}^N\int_{\Omega}c_i(x,t)dx}\\
\le C\sqrt{\mathcal D}\sqrt{\mathcal E + C}.
\ea
\la{rhsdis2}
\ee
We have thus
\be
\fr{d\mathcal F}{dt}\le C{\mathcal F +C},
\la{fineq}
\ee
and therefore 
\be
\sup_{0\le t\le T}\mathcal G(t)  \le \Gamma(T)
\la{Gbound}
\ee
and also
\be
\int_0^T\mathcal D_1(t)dt\le \Gamma(T)
\ee
for any $T>0$ with $\Gamma(T)$ depending only on $T$, initial data and boundary conditions. These estimates replace step 1, and the rest follows without new ideas like in the proof of existence for the uniform selective boundary conditions.
\beg{thm}\la{gensel} Under the assumptions of Theorem \ref{globlock} a unique global strong solution of the Nernst-Planck-Navier-Stokes system (\ref{pnpns}) with general selective boundary conditions
\be
\left\{
\ba
u_{\left |\right. \pa\Omega}=0,\\
\Phi_{\left |\right. \pa \Omega} = W(x),\\
c_i{\left |\right. S_i}= \gamma_i, \quad i=1, \dots, M,\\
(\pa_nc_i + z_ic_i\pa_n\Phi)_{\left |\right. \pa \Omega\setminus S_i} = 0,\quad i= 1,\dots, M,\\
(\pa_nc_i + z_ic_i\pa_n\Phi)_{\left |\right. \pa \Omega} = 0,\quad i= M+1,\dots , N
\ea
\right.
\la{lbc}
\ee
exists for any time $T$, and
\be
\sup_{0\le t\le T}[\|c_i(t)\|_{W^{2,q}(\Omega)} + \|\pa_tc_i(t)\|_{L^q(\Omega)}] \le C_q(T)
\la{maxlqTsel}
\ee
and
\be
\int_0^{T}[\|\pa_t u\|^2_{L^p(\Omega)} + \|u(t)\|_{W^{2,p}(\Omega)}^2]dt
\le U_p(T)
\la{maxregatsel}
\ee
hold. The constants $C_q(T)$ and $U_p(T)$ depend on the initial data and $T$.
\end{thm}

\section{Appendix A: Poisson-Boltzmann Equations}\la{poib}
We discuss here briefly the Poisson-Boltzmann equations encountered in the text. The subject is classical (\cite{friedman}, \cite{keller}).

We consider first the semilinear elliptic problem
\be
-\epsilon\D \Phi + G'(\Phi) = 0
\la{semiphi}
\ee
in the bounded domain $\Omega\subset\Rr^d$, with smooth boundary $\pa\Omega$.
The nonlinearity $G(\Phi)$ is given by
\be
G(\Phi) = \sum_{i=1}^N Z_i^{-1}e^{-z_i\Phi}
\la{gphi}
\ee
with $Z_i>0$ and $z_i\in\Rr$ given constants. The derivative $G^{'}(\Phi)$ is
\be
G^{'}(\Phi) = -\sum_{i=1}^N\fr{z_i}{Z_i}e^{-z_i\Phi}.
\la{gprimephi}
\ee
We note that $G$ is positive and convex.
The boundary conditions for $\Phi$ are (\ref{phibc}) with $W$ the boundary trace of a function $\widetilde W$,
\be
W = \widetilde{W}_{\left |\right. \pa\Omega}
\la{widetildeWag}
\ee
with $\widetilde{W}\in H^1(\Omega)\cap L^{\infty}(\Omega)$.  We let
\be
\mathcal A = \{ \Phi\in H^1(\Omega)\left |\right.\;  G(\Phi)\in L^{1}(\Omega), \; {\mbox{and}}\; \gamma_0(\Phi) = W\}
\la{adm}
\ee
where 
\be
\gamma_0(\Phi) =  \Phi_{\left |\right. \pa\Omega}
\la{trace}
\ee 
is the trace map $\gamma_0: H^1(\Omega) \to  H^{\fr{1}{2}}(\pa\Omega)$, 
and define, for $\Phi\in {\mathcal A}$,
\be
E(\Phi) = \int_{\Omega}\fr{\epsilon}{2}|\na\Phi|^2 + G(\Phi)dx.
\la{ephi}
\ee 
\beg{prop}\la{varprop} There exists $\Phi^*\in\mathcal A$ attaining the minimum of $E$:
\be
E(\Phi^*)= \min_{\Phi\in\mathcal A}E(\Phi).
\la{mine}
\ee
\end{prop}  
\noindent{\bf Proof}. Let $\alpha = \inf_{\Phi\in\mathcal A}E(\Phi)$. Because $E(\Phi)\ge 0$, there is no problem with the existence and finiteness of $\alpha\ge 0$. Let $\Phi_j\in{\mathcal A}$ be such that $\lim_{j\to\infty}E(\Phi_j) = \alpha$. The sequence $\Phi_j$ is bounded in $H^1(\Omega)$ and therefore the sequence $\Phi_j-{\widetilde W}$ is bounded in $H_0^1(\Omega)$. We can thus pass to a subsequence so that $\Phi_j -\widetilde W$ converge strongly in $L^2(\Omega)$, and consequently we can pass to a subsequence of $\Phi_j$ that converges weakly in $H^1(\Omega)$, strongly in $L^2(\Omega)$ and almost everywhere to a function $\Phi^*$. Because of the weak convergence in $H^1(\Omega)$ we have
\[
\int_{\Omega}|\na \Phi^*|^2 dx \le \liminf_{j\to\infty}\int_{\Omega}|\na \Phi_j|^2 dx.
\]
Because of the almost everywhere convergence and Fatou's lemma for the nonnegative functions $G(\Phi_j)$ we have that
\[
\int_{\Omega}G(\Phi^*)dx \le \liminf_{j\to\infty}\int_{\Omega}G(\Phi_j) dx,
\]
and because of the subadditivity of $\liminf$ we have
\[
E(\Phi^*) \le \lim\inf_{j\to\infty} E(\Phi_j) = \alpha.
\]
The inequalities above and the strong convergence in $L^2$ establish that $\Phi^*\in H^1(\Omega)$ and $G(\Phi^*)\in L^1(\Omega)$. Because the trace operator $\gamma_0$ is continuous between Hilbert spaces, hence weakly continuous, it follows that $\gamma_0(\Phi^*) = W$, and thus $\Phi^*\in {\mathcal A}$. This concludes the proof of the proposition.

We introduce
\be
\mathcal B = H_0^1(\Omega)\cap L^{\infty}(\Omega)
\la{bcal}
\ee
and observe that $\mathcal A + \mathcal B \subset \mathcal A$ (the sum of any element of $\mathcal A$ and any element of $\mathcal B$ belongs to $\mathcal A$). Then, fixing $\psi\in {\mathcal B}$ we observe that the function $s\mapsto E(\Phi^* + s\psi)$ is differentiable and  has a minimum at $s=0$. Carrying out the differentiation we arrive at the variational formulation:
\beg{prop} Let $\Phi^*$ be the minimum of $E$ on ${\mathcal A}$. Then, for any $\psi\in {\mathcal B}$ we have 
\be
\epsilon\int_{\Omega}\na \Phi^*\na \psi dx + \int_{\Omega}G'(\Phi^*)\psi dx = 0.
\la{varform} 
\ee
\end{prop}
We use now the variational formulation to gain regularity in a well established manner. We define
\be
\pa_h^i f(x) = \fr{f(x+he_i)-f(x)}{h}
\la{pah}
\ee
where $e_i =(0,\dots,0, 1, 0,\dots,0)$ is the canonical basis of $\Rr^d$, and $h\neq 0$. We note that
\be
(\pa_h^i)^* = -\pa_{-h}^i
\la{dual}
\ee
where the dual is with respect to the $L^2$ scalar product. We take a function
$\chi'_1$ of one variable that is smooth, even, compactly supported in the interval $[-2,2]$, is nonincreasing for positive $x$ and equals identically 1 on $[-1,1]$ and identically $0$ on $[-2,-\frac{3}{2}] \cup [\frac{3}{2}, 2]$. We define $\chi_1(x) = \int_0^x\chi_1'(s)ds$ and rescale $\chi'_M(x) = \chi'_1\left(\frac{x}{M}\right)$, and define  $\chi_M(x) = \int_0^x\chi'_M(s)ds$.  Note that $\chi_M$ is odd. We take another function $\chi\in C_0^{\infty}(\Omega)$. For any $1\le i\le d$, $M>1$ and $h\neq 0$, with $|h|< \fr{1}{2}{\mbox{dist}}({\mbox{supp}}\,\chi,\pa\Omega)$ we consider the test function
\be
\psi(x): = (\pa_h^i)^*[\chi(x)\chi_M(\pa_h^i\Phi^*(x))].
\la{test}
\ee
We easily check that $\psi\in {\mathcal B}$.  Now we apply the variational formulation (\ref{varform}). Let us describe the terms separately
\be
\ba
\epsilon\int_{\Omega}\na\Phi^*\na \psi dx 
= \epsilon\int_{\Omega}\pa_h^i\na\Phi^*(x) \na[\chi(x)\chi_M (\pa_h^i\Phi^*(x))]dx\\
= \epsilon\int_{\Omega}\na\pa_h^i\Phi^*(x) \chi(x) \na[\chi_M (\pa_h^i\Phi^*(x))]dx 
+ \epsilon\int_{\Omega}\na\pa_h^i\Phi^*(x)\na\chi(x) \chi_M (\pa_h^i \Phi^*(x))dx\\
= \epsilon\int_{\Omega}\na\pa_h^i\Phi^*(x)\chi(x)\chi_M' (\pa_h^i\Phi^*(x))\na\pa_h^i\Phi^*(x)dx 
+ \epsilon\int_{\Omega}\na\pa_h^i\Phi^*(x)\na\chi(x) \chi_M(\pa_h^i\Phi^*(x))dx\\
=\epsilon\int_{\Omega}|\na\pa_h^i\Phi^*(x)|^2\chi(x) \chi'_M(\pa_h^i\Phi^*(x))dx 
+ \epsilon\int_{\Omega}\na\chi(x) \na F_M(\pa_h^i\Phi^*(x))dx \\
= \epsilon\int_{\Omega}|\na\pa_h^i\Phi^*|^2\chi\chi'_M(\pa_h^i\Phi^*)dx 
- \epsilon\int_{\Omega}\D\chi  F_M(\pa_h^i\Phi^*)dx.
\ea
\la{nablas}
\ee
We used above the fact that $\pa_h^i$ and $\na$ commute. The function $F_M$ is given by
\be
F_M(\Phi) = \int_0^\Phi \chi_M(t)dt.
\la{fm}
\ee 
We note that, from our definitions
\be
F_M(\Phi)\le \fr{1}{2}\Phi^2.
\la{fmineq}
\ee
We obtained thus far:
\be
\epsilon\int_{\Omega}\na\Phi^*\na \psi dx 
= \epsilon\int_{\Omega}|\na\pa_h^i\Phi^*|^2\chi\chi'_M(\pa_h^i\Phi^*)dx 
- \epsilon\int_{\Omega}\D\chi  F_M(\pa_h^i\Phi^*)dx.
\la{nablaone}
\ee
Regarding the second term in (\ref{varform}) we have
\be
\int_{\Omega} G'(\Phi^*)\psi dx = \int_{\Omega}\pa_h^iG'(\Phi^*(x))\chi_M(\pa_h^i\Phi^*(x))\chi(x) dx.
\la{interg}
\ee
Now we observe that 
\be
\pa_h^iG'(\Phi^*(x)) = G''(S)\pa_h^i\Phi^*(x)
\la{s}
\ee
with $S$ some point on the segment $[\Phi^*(x), \Phi^*(x+h)]$. Observing that
\be
\Phi\chi_M(\Phi)\ge 0
\la{psim}
\ee
holds for any $\Phi$ we obtain from the convexity of $G$ that
\be
\int_{\Omega} G'(\Phi^*)\psi dx \ge 0.
\la{neg}
\ee
Adding (\ref{nablaone}) and (\ref{neg}), using (\ref{varform}) and (\ref{fmineq}) we obtain
\be
\epsilon\int_{\Omega}|\na\pa_h^i\Phi^*|^2\chi\chi'_M(\pa_h^i\Phi^*)dx \le \epsilon C_{\chi}\int_{\Omega}|\pa_h^i\Phi^*|^2dx.
\la{first}
\ee
Letting $M\to \infty$ and using the Lebesgue dominated convergence theorem we obtain
\be
\epsilon\int_{\Omega}|\na\pa_h^i\Phi^*|^2\chi dx 
\le \epsilon C_{\chi}\int_{\Omega}|\pa_h^i\Phi^*|^2dx\le C_{\chi} E(\Phi^*).
\la{second}
\ee
As a consequence, for any relatively compact subdomain $\Omega_1\subset\subset \Omega$ we have
\be
\epsilon \|\Phi^*\|_{H^2(\Omega_1)}^2 \le C E(\Phi^*).
\la{third}
\ee
This inequality implies, in $d=2, 3$, that $\Phi^*\in C^{\alpha}(\Omega_1)$. 
For higher dimensions we can show that $G'(\Phi^*)\in L^2_{loc}(\Omega)$.
In order to do so, we take the test function
\be
\psi(x) = \chi(x)\chi_M(G'(\Phi^*(x)))
\la{psinew}
\ee
with $\chi\in C_0^{\infty}(\Omega)$ and $\chi_M$ as above. It is easy to check that $\psi\in \mathcal B$ and thus we can apply (\ref{varform}).
We obtain
\be
\ba
0 = \epsilon\int_{\Omega}\na\Phi^*\na\psi + \int_{\Omega}G'(\Phi^*)\psi 
= \epsilon\int_{\Omega} |\na\Phi^*|^2\chi_M'(G'(\Phi^*)) G''(\Phi^*)\chi dx \\
+ \epsilon\int_{\Omega} \na\Phi^* \chi_M(G'(\Phi^*))\na\chi dx + \int_{\Omega} G'(\Phi^*)\chi_M(G'(\Phi^*))\chi dx.
\ea
\la{intergm}
\ee
Now we note that
\be
x\chi_M(x) \ge \fr{1}{2}\chi_M^2(x)
\la{chimeneq}
\ee
which can be verified easily by differentiation, noticing that, in view of the fact that the functions are even it is enough to check for nonnegative $x$, and using the fact that $\chi_M(x) \le x$ for nonnegative $x$.
We obtain, using a Schwartz inequality:
\be
 \epsilon\int_{\Omega} |\na\Phi^*|^2\chi_M'(G'(\Phi^*)) G''(\Phi^*)\chi dx 
+ \frac{1}{2}\int_{\Omega} G'(\Phi^*) \chi_M(G'(\Phi^*))\chi dx \le
C\epsilon^2\int_{\Omega}\frac{|\na\chi|^2}{\chi}|\na\Phi^*|^2dx.
\la{intergm2}
\ee
Letting $M\to\infty$ and using the fact that $x\chi_M(x)$ is a nonnegative function which is nondecreasing in $M$, we obtain from the monotone convergence theorem and the convexity of $G$
\be
\int_{\Omega}G'(\Phi^*)^2\chi dx \le C\epsilon^2\int_{\Omega}\frac{|\na\chi|^2}{\chi}|\na\Phi^*|^2dx.
\la{gprimesquarechi}
\ee
Thus, $G'(\Phi^*)\in L^2_{loc}(\Omega)$ and, using elliptic regularity we can bootstrap and obtain bounds for higher derivatives in any dimension $d$. We will not pursue this here.

\beg{thm}\la{phistarnthm} 
Let $\Omega\subset \Rr^d$, $d=2,3$ be bounded domain with smooth boundary and let $W$ be a smooth enough function on $\pa \Omega$ (for instance $W\in H^s(\pa\Omega)$ with $s\ge\fr{3}{2}$). Then there exists a unique weak solution 
$\Phi^*\in\mathcal A\cap H^2_{loc}(\Omega)$ of 
\be
-\epsilon \D\Phi^* = -G'(\Phi^*)
\la{phistarneq}
\ee
with boundary condition (\ref{phistarbc}), with
\be
G(\Phi^*) = \sum_{i=1}^N Z_i^{-1}e^{-z_i\Phi^*}
\la{Gn}
\ee
and with given positive numbers $Z_i>0$ in the sense of (\ref{varform}). If $W\in W^{\fr{3}{2}, p}(\pa\Omega)$ with $p>d$ then $\Phi^*\in W^{2,p}(\Omega)$ and consequently $\Phi^*\in W^{1,\infty}(\Omega)$.
\end{thm}
\noindent{\bf Proof.} The existence and interior regularity have been established. The formal calculation for the uniqueness is simple: If $\Phi^*_i\in\mathcal A$, for $i=1,2$, are two weak solutions then
\be
-\epsilon \D(\Phi^*_1-\Phi^*_2) + G'(\Phi^*_1)-G'(\Phi^*_2)=0
\la{diffeq}
\ee
with $\Phi^*_1-\Phi^*_2 \in H_0^1(\Omega)$. Taking the scalar product with $\Phi^*_1-\Phi^*_2$ and observing that 
\be
(G'(\Phi^*_1)-G'(\Phi^*_2))(\Phi^*_1-\Phi^*_2)\ge 0
\la{conv}
\ee
holds pointwise because of the convexity of $G$, we obtain that 
\be
\int_{\Omega}|\na(\Phi^*_1-\Phi^*_2)|^2dx = 0.
\la{diffzero}
\ee
The rigorous argument is as follows: by the interior regularity of solutions,
(\ref{diffeq}) holds almost everywhere in $\Omega$, and the inequality (\ref{conv}) is pointwise true. Therefore the function
\[
(\D(\Phi^*_1-\Phi^*_2))(\Phi^*_1-\Phi^*_2)
\]
which a priori is known to be in $L^1_{loc}(\Omega)$ is nonnegative almost everywhere. Thus, from interior regularity, denoting $\psi = \Phi^*_1-\Phi^*_2$ we have
\[
|\na \psi|^2 \le \fr{1}{2}\D\psi^2
\]
almost everywhere. The left hand side is in $L^1(\Omega)$, as $\psi\in H_0^1(\Omega)$, and the right hand side is in $L^{1}_{loc}(\Omega)$ by interior regularity. Taking now $w_1$, the positive eigenfunction corresponding to the first eigenvalue of $-\D$ with homogeneous boundary conditions, we obtain
\[
\int_{\Omega}w_1|\na \psi|^2 dx \le \fr{1}{2}\int_{\Omega}w_1\D\psi^2dx = -\fr{\lambda_1}{2}\int_{\Omega}w_1\psi^2dx.
\]
The integration by parts is allowed because $\psi\in H_0^1(\Omega)$ and $w_1$ can be approximated in $H^1(\Omega)$ by $C_0^{\infty}$ functions. This shows that $\psi=0$, because, as it is well known, $w_1(x)\ge Cd(x)>0$ where $d(x)$ is the distance from $x$ to the boundary of the domain.

Let us make a few remarks about (\ref{phistarneq}). If $W_1(x)\le W_2(x)$ are two boundary conditions, and if $\Phi^*_1, \Phi^*_2$ denote the corresponding solutions, (assumed to be continuous up to the boundary) it follows from the maximum principle that $\Phi_1^*(x)\le \Phi_2^*(x)$ everywhere. Indeed, from $G''>0$  it follows from the equations that $\Phi_1-\Phi_2$ cannot attain its maximum in the interior of the domain. 

If $Z_i$ together with $z_i$ satisfy the neutrality condition
\be
\sum_{i=1}^N \fr{z_i}{Z_i} = 0,
\la{neutn}
\ee
then $G'(0)=0$. In this case, using the fact that $G'(0)=0$ and the fact that $0$ solves the equation with zero boundary conditions, it follows that if $W(x)\ge 0$ on the boundary, the corresponding solution is nonnegative $\Phi^*(x)\ge 0$. Then, considering $M= \max|W(x)|$ on the boundary it follows that $\Phi^*(x)\le \Phi_M^*(x)$ where $\Phi^*_M$ solves the problem (\ref{phistarneq}) with constant boundary condition equal to $M$. Because $\Phi^*_M(x)\ge 0$
and $G'(\Phi)>0$ for $\Phi>0$ it follows again from the maximum principle that
$\Phi_M^*(x)\le M$. Therefore, for any $W$ we have 
\be
-M\le\Phi^*(x)\le M.
\la{phistarm}
\ee  
This bound is remarkable in that it does not depend on $z_i, Z_i$, once the neutrality condition is assumed. The considerations above can be made rigorous, for instance by adding a small multiple of $ G^2(\Phi)$ to the variational problem, and then removing it. The minimization of
\be
\int_{\Omega}\left[\fr{\epsilon}{2}|\na\Phi|^2 + G_r(\phi)\right]dx
\la{endelta}
\ee
with
\be
G_r(\Phi) = G(\Phi) + r G^2(\Phi)
\la{Gr}
\ee
with $r>0$ on the corresponding admissible set ${\mathcal A}_r = \{\Phi\in H^1(\Omega)\left|\right. G_r(\Phi)\in L^1(\Omega), \gamma_0(\Phi)= W\}$ yields bounded solutions with the same $L^{\infty}$ bounds, and their regularity up to the boundary is classical. Removing $r$ we deduce the bounds (\ref{phistarm}) for $\Phi^*$ and then again we can apply classical results to obtain regularity up to the boundary.

Let us provide here an explicit calculation for a one dimensional case, similar to  to one used in \cite{howard} in a half-space, using the neutrality condition. Let 
\be
-\epsilon \Phi'' + G'(\Phi) = 0
\la{onedeq}
\ee
on the interval $[0,H]$
with boundary conditions
\be
\Phi (H) = W, \quad \Phi(0) = 0
\la{onedbc}
\ee
with $W>0$.  Multiplying (\ref{onedeq}) by $\Phi'$ and integrating once we obtain
\be
\epsilon (\Phi')^2 = 2(G(\Phi) - A)
\la{square}
\ee
with $A$ a constant of integration. If we are to have smooth solutions, $A$ must not exceed the minimum of $G(\Phi)$ on the interval. Now $G$ is convex and the global minimum of $G$ is $G(0)$ because $G'(0)= 0$. Because $0$ is in the range of $\Phi$ (it is a boundary condition) it follows that the minimum of $G(\Phi)$ is $G(0)$. We write  $A= G(0) -\alpha^2$ with $\alpha\ge 0$. We choose $\alpha$ such that
\be
\int_0^W\fr{d\Phi}{\sqrt{G(\Phi) - G(0) + \alpha^2}}d\Phi = \sqrt{\fr{2}{\epsilon}}H.
\la{alphaeq}
\ee
The fact that we can solve this equation requires a small argument, based on the fact that when $\alpha=0$ the integral diverges and the fact that $G$ is convex. Thus
\[
C_1\Phi^2 \le G(\Phi) - G(0) 
\]
for $\Phi\in [0, W]$ because of convexity, and
\[
G(\Phi) - G(0)\le C_2\Phi^2
\]
for $\Phi\in [0, \Phi_0]$ because of continuity of the second derivative of $G$, with $C_1>0$, $C_2$ and $\Phi_0$ independent of $W$. Therefore part of the integral in (\ref{alphaeq}) is bounded below by  
\[
\int_0^{\Phi_0}\fr{d\Phi}{\sqrt{C_2\Phi^2 + \alpha^2}} \ge \fr{1}{\sqrt{2C_2}}\log\left(\fr{\Phi_0\sqrt{C_2}}{\alpha}\right)
\]
and the rest from above by
\[
\int_{\Phi_0}^W\fr{d\Phi}{\sqrt{C_1\Phi^2 + \alpha^2}} \le\fr{1}{\sqrt{C_1}}\ \log\left(\fr{W}{\Phi_0}\right).
\]
The sum therefore can be made arbitrarily large, as $W$ is fixed (even if it depends on $\epsilon$) and $\alpha$ is chosen small enough. On the other hand, if $\alpha$ is large enough, then the integral on the left hand side of (\ref{alphaeq}) can be made arbitrarily small. Thus, as $\alpha$ is varied, the range of the integral contains the target value in the right hand side of (\ref{alphaeq}). 

We then set
\be
P(\Phi) = \int_0^\Phi\fr{d\Psi}{\sqrt{G(\Psi) - G(0) + \alpha^2}}d\Psi 
\la{psi}
\ee
and 
\be
\Phi^* (y) = P^{-1}\left(\sqrt{\fr{2}{\epsilon}} y\right)
\la{onedphistar}
\ee
and conclude the construction.

Let us turn now to the equation (\ref{phinftyeq}) which is the Poisson-Boltzmann equation for the case of blocking boundary conditions for the ionic species,
namely,
\be
-\epsilon\D\Phi^*= \sum_{i=1}^N z_i I_i^0 \fr{e^{-z_i\Phi^*}}{\int_{\Omega} e^{-z_i\Phi^*}}.
\la{phistarbl}
\ee
The constants $I_i^0$ are given positive numbers and the  boundary conditions are  (\ref{phibc}).
This is special case of (\ref{phistarneq}), (\ref{Gn}), in which
\be
Z_i = (I_i^0)^{-1}\int_{\Omega}e^{-z_i\Phi^*}dx.
\la{ziphinftyn}
\ee
Regularity of weak solutions follows from the fact that the equation is semilinear elliptic and the smoothness of the boundary and of the boundary conditions. Uniqueness was shown in Theorem \ref{uniqphinfty} and existence is a consequence of Theorem \ref{asy}. There are several other approaches to show existence. Showing that the equation (\ref{ziphinftyn}) can be solved after finding $\Phi^*(Z_1, \dots Z_N)$ solutions of (\ref{phinftyeq}) is a nontrivial possible route.  A  proof of existence using the fact that solutions are critical points of the energy
\be
\fr{\epsilon}{2}\int_{\Omega} |\na \Phi(x)|^2dx + \sum_{i=1}^NI_i^{0}\log\left (\int_{\Omega}e^{-z_i\Phi(x)}dx\right)
\la{nonlocale}
\ee
is also possible. This is the approach in \cite{friedman}  who solve a special case (albeit with slightly different boundary conditions). The energy is bounded below by Jensen's inequality and an approximation is used to control the exponential integrals in the logarithms. 

Finally, we also consider the Poisson - Boltzmann equation for the uniform selective boundary conditions,
\be
-\epsilon\D\Phi^* = \sum_{i=1}^M z_iZ_i^{-1}e^{-z_i\Phi^*} + \sum_{i=M+1}^N z_iI_i^0
\fr{e^{-z_i\Phi^*}}{\int_{\Omega} e^{-z_i\Phi^*}dx}
\la{mixedeqphistar}
\ee
with boundary conditions (\ref{phibc}). The existence of solutions follows from Theorem \ref{uniblock}, regularity follows from the semilinear elliptic character and  the uniqueness follows in the manner of Theorem \ref{uniqphinfty}. A direct existence proof can be constructed using the fact that solutions are critical points of
\be
\fr{\epsilon}{2}\int_{\Omega} |\na \Phi(x)|^2dx + \int_{\Omega}\sum_{i=1}^MZ_i^{-1}e^{-z_i\Phi(x)}dx + \sum_{i=M+1}^NI_i^{0}\log\left (\int_{\Omega}e^{-z_i\Phi(x)}dx\right)   
\la{mixe}
\ee
with $Z_i>0$ and $I_i^0>0$ given numbers. In fact (\ref{mixedeqphistar}) include both (\ref{phistarneq}) (\ref{Gn}), when $M=N$  and (\ref{phistarbl}) when $M=0$.  

It is interesting to note that the linearization of equation (\ref{mixedeqphistar}) at a state $\Phi$ is the linear elliptic nonlocal operator
\be
L_\Phi(\psi) = -\epsilon \D\psi + G''(\Phi)\psi + \sum_{i=M+1}^Nz_i^2I_i^0(\psi -(\psi,p_i)_{L^2(\Omega)})p_i
\la{lpsi}
\ee
where $G(\Phi) = \sum_{i=1}^M Z_i^{-1}e^{-z_i\Phi}$ and  $p_i = \fr{e^{-z_i\Phi(x)}}{\int_{\Omega}e^{-z_i\Phi(x)}dx}$. This operator with domain $H^2(\Omega)\cap H_0^1(\Omega)$ is selfadjoint in $L^2(\Omega)$, positive and invertible when $\Phi\in L^{\infty}(\Omega)$. These properties can be used to produce a nontrivial Newton iteration procedure for computing solutions of (\ref{mixedeqphistar}).
\section{Appendix B}\la{locex}
We sketch here for the sake of completeness our proof of Theorem \ref{loc}. Local existence based on methods of maximal regularity was presented in (\cite{bothe}.)

We consider an iteration:
\be
\pa_t c_i = D_i(\D c_i + z_i\div(c_i\na \Phi_o)) - u\cdot\na c_i
\la{ciit}
\ee
with 
\be
-\epsilon\D \Phi_o = \rho_o
\la{phiit}
\ee
and
\be
\pa_t u + u\cdot\na u +\na p = \nu\D u - (k_B\T) \rho_o\na\Phi_o, \quad \div u = 0,
\la{nsit}
\ee
boundary conditions
\be
(\na c_i + z_ic_i\na \Phi_o)_{\left |\right. \; \pa\Omega}\;  \cdot n = 0
\la{ciitbc} 
\ee
\be
{\Phi_o}_{\left |\right. \; \pa\Omega} = W
\la{phiitbc}
\ee
and 
\be
u_{\left |\right. \; \pa\Omega} =0.
\la{uitbc}
\ee
We are assuming that $\rho_o(x,t)$ is given by a previous calculation, and we are interested in inductive bounds. We do not mention explicitly the counting index of the iteration. We observe that the linear equations (\ref{ciit}) with time dependent boundary conditions (\ref{ciitbc}) are equivalent to the linear equations
\be
\pa_t{\widetilde{c_i}} = D_i\D {\widetilde{c_i}} - (u + D_iz_i\na\Phi_o)\na {\widetilde{c_i}} + z_i((\pa_t + u\cdot\na)\Phi_o){\widetilde{c_i}} 
\la{Cieq}
\ee
with homogeneous Neumann boundary conditions
\be
\pa_n {{\widetilde{c_i}}}_{\left |\right. \; \pa\Omega} = 0
\la{Cibc}
\ee
for the dependent variable
\be
{\widetilde{c_i}} = c_ie^{z_i\Phi_o}.
\la {Ci}
\ee 
This observation clarifies the nature of the equations: Obviously, if $\rho_o$, and consequently $\Phi_o$, $u$ are smooth, then ${\widetilde{c_i}}$, and consequently $c_i$ are smooth.  This allows us to perform calculations on the preferred form (\ref{ciit}). 
We start by estimating norms $\|c_i(t)\|_{L^p(\Omega)}$ for $p>d$. Integrating by parts and using the boundary conditions we have
\be
\fr{1}{p(p-1)D_i}\fr{d}{dt}\|c_i(t)\|_{L^p(\Omega)}^p = -\int_{\Omega}(|\na c_i|^2 + z_ic_i\na\Phi_o\cdot\na c_i)c_i^{p-2}dx
\la{ciitlpeq}
\ee
and therefore
\be
\fr{1}{p(p-1)D_i}\fr{d}{dt}\|c_i(t)\|_{L^p(\Omega)}^p + \fr{1}{2}\int_{\Omega}
|\na c_i|^2c_i^{p-2}dx \le \fr{z_i^2}{2}\|\na \Phi_o(t)\|_{L^{\infty}(\Omega)}^2\|c_i(t)\|_{L^p(\Omega)}^p .
\la{ciitlpineq}
\ee
Consequently, we have
\be
\|c_i(t)\|_{L^p(\Omega)}\le e\|c_i(0)\|_{L^p(\Omega)}
\la{ecipl}
\ee
for all $0\le t\le T$, provided
\be
\int_0^T\|\na \Phi_o(t)\|_{L^{\infty}(\Omega)}^2dt \le \fr{2}{(p-1)D_iz_i^2}.
\la{cond1}
\ee
Let us consider the inductive situation, when
\be
\rho_o(x,t) = \sum_{i=1}^N z_i c_i^{o}(x,t),
\la{rhooco}
\ee
and let us assume the time interval $[0, T_0]$ we have that 
\be
\sup_{0\le t\le T_0}\|c_i^{o}(t)\|_{L^p(\Omega)} \le A_p
\la{ap}
\ee
for some $p>d$. Using elliptic regularity, we have that
\be
\|\na \Phi_o\|_{L^{\infty}(\Omega)}\le C_{\Omega}[\|\rho_o\|_{L^p(\Omega)} + \|W\|_{W^{\fr{3}{2}, p}(\pa\Omega)}]
\la{regphio}
\ee
and, taking into account (\ref{rhooco}) and the assumption (\ref{ap}) we have that
\be
\sup_{0\le t\le T_0}\|\na \Phi_o\|_{L^{\infty}(\Omega)}\le C_W(A_p +1) 
\la{naphioind}
\ee
where we took
\be
C_W = C_{\Omega}(\sum_{i=1}^N|z_i| +  \|W\|_{W^{\fr{3}{2}, p}(\pa\Omega)}),
\la{CW}
\ee
a constant that depends only on the data of the problem. The condition (\ref{cond1}) is then satisfied if
\be
TC_W^2(A_p+1)^2 \le  \fr{2}{(p-1)D_iz_i^2},
\la{cond2}
\ee
and, if that is the case, we guarantee (\ref{ecipl}) on the interval $[0,T]$. Therefore, choosing
\be
A_p = e\|c_i(0)\|_{L^p(\Omega)}
\la{apchoice}
\ee
we conclude that the assumption (\ref{ap}) is preserved in the iteration,
\be
\sup_{0\le t\le T_0}\|c_i(t)\|_{L^p(\Omega)}\le e\|c_i(0)\|_{L^p(\Omega)} = A_p,
\la{cilpb}
\ee
if 
\be
T_0 \le \fr{2}{(p-1)\max_i(D_i z_i^2)}C_W^{-2}\left(e\|c_i(0)\|_{L^p(\Omega)} + 1\right)^{-2}.
\la{tzero}
\ee
 Let us note that from (\ref{ciitlpineq}) we have also
\be
\int_0^{T_0} \int_{\Omega}|\na c_i|^2c_i^{p-2}dxdt \le \fr{2p+1}{D_ip(p-1)}A_p^p.
\la{gradcilpb}
\ee
In order to provide further inductive information we require that $p>2d$ and that $T_0$ satisfies the constraint (\ref{tzero}) with a possibly larger constant $M$,
\be
T_0\le M^{-1}\left(e\|c_i(0)\|_{L^p(\Omega)} + 1\right)^{-2}.
\la{tzerom}
\ee
At this point we have required only
\be
M\ge  \fr{(p-1)\max_i(D_i z_i^2)C_W^{2}}{2}.
\la{mone}
\ee
We provide below the justification for the additional requirement
\be
M\ge 16e^2(q-1) \max_i(D_iz_i^2)C_z^2|\Omega|^{\fr{2(p-q)}{pq}}.
\la{mtwo}
\ee

We remark that the condition (\ref{tzerom}) depends only on the norms $\|c_i(0)\|_{L^p(\Omega)}$ of the initial data and on the parameters of the problem, but not on the iteration step, nor on higher regularity data, or velocity initial data.

The equation (\ref{ciit}) can be written as
\be
\pa_t c_i + \div j_i = 0
\la{ciitdivj}
\ee
with 
\be
j_i = -D_i(\na c_i + z_ic_i\na\Phi_o) + u c_i.
\la{jiit}
\ee
We take the time derivative and use the fact that the boundary conditions imply 
\be
{j_i}_{\left |\right.\;\pa \Omega}\;\cdot n = 0.
\la{jitbc}
\ee
The time derivative $\pa_t c_i $ obeys thus
\be
\pa_t(\pa_t c_i) + \div \pa_t(j_i) = 0
\la{patcieq}
\ee
with boundary condition
\be
{\pa_tj_i}_{\left |\right.\;\pa \Omega}\;\cdot n = 0.
\la{jtbc}
\ee
We multiply (\ref{patcieq}) by $(\pa_t c_i)|\pa_t c_i|^{q-2}$ for some $q\ge 2$ and integrate by parts. We obtain
\be
\fr{1}{q}\fr{d}{dt}\int_{\Omega} |\pa_t c_i|^q dx = \int_{\Omega}(\pa_t j_i \cdot\na ((\pa_t c_i)|\pa_t c_i|^{q-2})dx.
\la{patcilpeq}
\ee 
This yields
\be
\ba
\fr{1}{q(q-1)}\fr{d}{dt}\int_{\Omega} |\pa_t c_i|^q dx = -D_i\int_{\Omega}\left[|\na\pa_t c_i|^2 +z_i(\pa_t c_i) \na\Phi_o\cdot\na \pa_t c_i\right]|\pa_t c_i|^{q-2}dx\\
+\int_{\Omega}c_i\left[(\pa_t u - D_iz_i\na\pa_t\Phi_o)\cdot \na\pa_t c_i\right]|\pa_t c_i|^{q-2}dx.
\ea
 \la{patcilpeq1}
\ee
Consequently we have
\be
\ba
\fr{1}{q(q-1)}\fr{d}{dt}\|\pa_t c_i\|_{L^q(\Omega)}^q  +\fr{D_i}{4}\int_{\Omega}|\na\pa_t c_i|^2 ]|\pa_t c_i|^{q-2}dx\\
\le \fr{D_iz_i^2}{2}\|\na\Phi_o(t)\|_{L^{\infty}}^2\|\pa_t c_i\|_{L^q(\Omega)}^q + \fr{1}{D_i}\|c_i(\pa_t u - D_iz_i\na\pa_t\Phi_o)\|_{L^q(\Omega)}^2\|\pa_tc_i\|_{L^q(\Omega)}^{q-2},
\ea
\la{patcilpin}
\ee
where we used a H\"{older} inequality with exponents $2, q, \fr{2q}{q-2}$ and Schwartz inequalities. We have from (\ref{patcilpin})
\be
\fr{d}{dt}\|\pa_t c_i\|_{L^q(\Omega)}^2\le (q-1)D_iz_i^2\|\na\Phi_o(t)\|_{L^{\infty}}^2 \|\pa_t c_i\|_{L^q(\Omega)}^2 + \fr{2(q-1)}{D_i}\|c_i(\pa_t u - D_iz_i\na\pa_t\Phi_o)\|_{L^q(\Omega)}^2.
\la{patcilpine}
\ee
From (\ref{cond1}) we obtain that
\be
\|\pa_t c_i(t)\|_{L^q(\Omega)}^2 \le e^2\|\pa_t c_i(0)\|_{L^q(\Omega)}^2 +
\fr{4e^2(q-1)}{D_i}\int_0^{T_0}\left[\|c_i\pa_t u\|_{L^q(\Omega)}^2 + D_i^2z_i^2\|c_i\na\pa_t\Phi_o)\|_{L^q(\Omega)}^2\right]dt
\la{patcilpind}
\ee
holds for all $t\le T_{0}$. We treat the two integral terms in the right hand side of (\ref{patcilpind}) differently. Because
\be
-\epsilon\D\pa_t\Phi_o = \pa_t\rho_o
\la{patpoi}
\ee
with boundary condition
\be
\pa_t{\Phi_o}_{\left | \right. \;\pa\Omega} = 0
\la{patphibc}
\ee
we have, from elliptic regularity
\be
\|\pa_t\Phi_o(t)\|_{W^{1,\infty}(\Omega)}\le C_{\Omega}\|\pa_t\rho_o\|_{L^q(\Omega)}.
\la{patel}
\ee
Let us assume that
\be
\sup_{0\le t\le T_0}\|\pa_t c_i^o(t)\|_{L^q(\Omega)} \le B_q.
\la{patcilqind}
\ee
Then it follows that
\be
\sup_{0\le t\le T_0}\|\pa_t\Phi_o(t)\|_{W^{1,\infty}(\Omega)}\le C_{z}B_q
\la{napatphi}
\ee
with
\be
C_{z} = C_{\Omega}\sum_{i=1}^N |z_i|. 
\la{cz}
\ee
Thus for the second integral term in the right hand side of (\ref{patcilpind}) we obtain that
\be
4e^2(q-1)\int_0^{T_0} D_iz_i^2\|c_i\na\pa_t\Phi_o)\|_{L^q(\Omega)}^2dt \le
4e^2(q-1) D_iz_i^2C_z^2B_q^2A_q^2T_0
\la{intonebound}
\ee
which implies
\be
4e^2(q-1)\int_0^{T_0} D_iz_i^2\|c_i\na\pa_t\Phi_o)\|_{L^q(\Omega)}^2dt \le \fr{1}{4}B_q^2
\la{first quarter}
\ee
if $q\le p$, in view of (\ref{mtwo} ) of condition (\ref{tzerom}).
For the first integral term we use $p=2q$ and bound
\be
\fr{4e^2(q-1)}{D_i}\int_0^{T_0}\|c_i\pa_t u\|_{L^q(\Omega)}^2 dt\le \fr{4e^2(q-1)}{D_i}A_p^2\int_0^{T_0}\|\pa_t u\|^2_{L^p(\Omega)}dt .
\la{intwo}
\ee
Now we use the bound (\cite{solo}, \cite{solo1})
\be
\int_0^{T_0}\|\pa_t u\|^2_{L^p(\Omega)}dt\le C\left (\|u(0)\|_{W^{1,p}(\Omega)}^2 + \int_0^{T_0}\|F(t)\|_{L^p(\Omega)}^2dt\right)
\la{maxreg}
\ee
which is valid on any time interval in $d=2$ and on a short time interval, independent of iteration in $d=3$. Here 
\be
F = -(k_B\T) \rho_o \na\Phi_o
\la{fnso}
\ee
obeys in view of (\ref{naphioind}), (\ref{cilpb})
\be
\|F(t)\|_{L^p(\Omega)} \le (k_B\T)C_zA_pC_W(A_p+1)
\la{fbound}
\ee
and consequently 
\begin{align}
&\fr{4e^2(q-1)}{D_i}\int_0^{T_0}\|c_i\pa_t u\|_{L^q(\Omega)}^2 dt\\
&\le \fr{4e^2(q-1)}{D_i}A_p^2 C\left[ \|\pa_t u(0)\|_{L^p(\Omega)}^2 +  (k_B\T)^2C_z^2A_p^2C_W^2(A_p+1)^2T_0\right].
\la{intwob}
\end{align}
Consequently, using (\ref{tzerom}) 
\be
\fr{4e^2(q-1)}{D_i}\int_0^{T_0}\|c_i(\pa_t u)\|_{L^q(\Omega)}^2 dt\le \fr{1}{4}B_q^2
\la{intowbou}
\ee
if we impose
\be
B_q^2\ge C_1A_p^2(\|u(0)\|_{W^{1,p}(\Omega)}^2 + A_p^2) + 2e^2 \|\pa_t c_i(0)\|_{L^q(\Omega)}^2
\la{bqreq}
\ee
with $C_1$ depending only on the parameters of the problem. Then, returning
to (\ref{patcilpind}) we have
\be
\|\pa_t c_i(t)\|_{L^q}^2 \le B_q^2
\la{patcilqbound}
\ee
for all $t\le T_0$.
We return now to the equation (\ref{ciit}) written as
\be
-\D c_i = -\fr{1}{D_i}\pa_t c_i  + (z_i\na\Phi_o -\fr{1}{D_i} u)\na c_i -\fr{z_i}{\epsilon}\rho_o c_i
\la{ciitlap}
\ee
and esitimate the right hand side in $L^q$ using (\ref{rhooco}), (\ref{ap}), (\ref{naphioind}), (\ref{cilpb}) and (\ref{patcilpind})
\be
\|\D c_i(t)\|_{L^q(\Omega)}\le \fr{1}{D_i}B_q + \left(|z_i|C_W(A_p+1) + \fr{1}{D_i}\|u(t)\|_{L^{\infty}(\Omega)}\right)\|\na c_i(t)\|_{L^q(\Omega)} + \left(\fr{|z_i|}{\epsilon}\sum_{j=1}^N|z_j|\right) A_p^2.
\la{delindint}
\ee
In order to finish we use  the variables $\widetilde{c}_i$  defined in (\ref{Ci}) which obey homogeneous Neumann boundary conditions. They obey therefore elliptic bounds
\be
\|{\widetilde{c}}_i\|_{W^{2,q}(\Omega)}\le C_{\Omega}\left( \|\D\widetilde c_i\|_{L^q(\Omega)} + \|\widetilde c_i\|_{L^q(\Omega)}\right)
\la{ellwidetc}
\ee
and, integrating by parts we see that
\be
\|\na\widetilde c_i\|_{L^q(\Omega)}^2 \le C_{\Omega}\|\widetilde{c_i}\|_{L^q(\Omega)}\left(\|\D\widetilde{c}_i\|_{L^q(\Omega)}  + \|\widetilde c_i\|_{L^q(\Omega)}\right).
\la{sqrtilde}
\ee
Returning to the variables $c_i$ we have, in view of (\ref{ap}) and (\ref{naphioind})
\be
\|\widetilde c_i\|_{L^q(\Omega)}\le e^{|z_i|C_W(A_p+1)}\|c_i\|_{L^q(\Omega)}
\le  e^{|z_i|C_W(A_p+1)}A_p|\Omega|^{\fr{1}{q}-\fr{1}{p}} = H_0(A_p),
\ee
and similarly,
\be
\|\na c_i\|_{L^q(\Omega)} \le H_1(A_p)\left(\|\na \widetilde{c}_i\|_{L^q(\Omega)} + \|c_i\|_{L^q(\Omega)}\right)
\la{grads}
\ee
and
\be
\|\D {\widetilde {c}}_i\|_{L^q(\Omega)}\le H_2(A_p)(\|\D c_i\|_{L^q} + 1)
\la{hap}
\ee
with $H_1$ and $H_2$ explicit functions of $A_p$. Therefore, from 
(\ref{delindint}) we obtain
\be
\sup_{0\le t\le T_0}\|c_i\|_{W^{2,q}(\Omega)} \le H_3(A_p, B_q, \|u_0\|_{W^{1,q}(\Omega)})
\la{fin}
\ee
where $H_3$ is an explicit positive continuous function, nondecreasing in each of its arguments, and depending also on the parameters $z_i, \nu, \epsilon$ but not on the iteration step. 

We construct thus by induction a sequence of solutions of linear equations (\ref{ciit}), (\ref{phiit}), (\ref{nsit}) which obey uniform bounds (\ref{cilpb}) on a common interval of time $[0, T_0]$, determined by the condition (\ref{tzerom}) with (\ref{mone}) and (\ref{mtwo}). We have also the bounds for higher derivatives (\ref{patcilqbound}), (\ref{maxreg}), (\ref{fin}). Passing to the limit in the sequence is straightforward and yields a short time solution with the stated bounds.

\section{Conclusion}
We proved global existence of solutions for two dimensional Nernst-Planck-Navier-Stokes equations in bounded domains for arbitrary large initial data, arbitrary valences, voltages, different species diffusivities, any dielectric constant and arbitrary Reynolds numbers, in the cases of both blocking and general selective boundary conditions. Convergence to uniquely determined Boltzmann states and zero fluid velocity occurs not only for blocking boundary conditions, but also for uniform selective conditions. The latter include complex nontrivial configurations in which large voltage differences can be applied.  

\vspace{.5cm}

\noindent{\bf{Acknowledgment.}} The work of PC was partially supported by NSF grant DMS-1713985. 


\begin{thebibliography}{99}
\bibitem{biler} P. Biler, J. Dolbeault. Long time behavior of solutions to Nernst-Planck and Debye-Hckel drift-diffusion systems. Ann. Henri Poincare {\bf{1}}(2000) 461-472.
\bibitem{bothe} D. Bothe, A. Fischer, J. Saal, Global well-posedness and stability of electrokinetic flows, SIAM J. Math. Anal, {\bf 46} 2, (2014), 1263-1316.
\bibitem{choi} Y.S. Choi, and R. Lui, Multi-Dimensional Electrochemistry Model, Arch. Rational Mech. Anal. {\bf{130}} (1995), 315-342.
\bibitem{davidson} S. M. Davidson, M. Wessling, A. Mani, On the dynamical regimes of pattern-accelerated electroconvection, Scientific Reports {\bf 6} 22505 (2016) doi:19.1039/srep22505
\bibitem {fischer} A. Fischer, J. Saal, Global weak solutions in three space dimensions for electrokinetic flow processes, Journal of Evolution Equations, {\bf 17} 1, (2017), 309-333.
\bibitem{friedman} A. Friedman, K. Tintarev, Boundary asymptotics for solutions of the Poisson-Boltzmann equation, J. Diff. Eqn {\bf{69}} (1987), 15-38.
\bibitem{gajewski} H. Gajewski, K. Groger, Reaction-diffusion processes of electrically charged species. Math. Nachr., {\bf{177}} (1996), 109-130.
\bibitem{solo}G. Grubb, V. Solonnikov, Boundary value problems for nonstationary Navier-Stokes equations treated by pseudo-differential methods, Math. Scand. {\bf{69}} 217-290, (1991).
\bibitem{howard} A. Gupta, H. Stone, Consequences of asymmetry in electrolyte valence on diffuse charge dynamics, preprint 2018.
\bibitem{jerome} J. W. Jerome, Analytical approaches to charge transport in a moving medium, Transport Theory Stat Phys. {\bf 31} (2002) 333-366.
\bibitem{jeromesacco} J. W. Jerome, R. Sacco, Global weak solutions for an incompressible charged fluid with multi-scale couplings: Initial-boundary-value problem, Nonlinear Analysis, Theory, Models and Applications, {\bf 71} 12 (2009),
e2487-e2497.
\bibitem{keller}J. B. Keller, Electrohydrodynamics I, The equilibrium of a charged gas in a container, J. Rational Mech. Anal. {\bf{5}} (1956), 715-724.
\bibitem{rubibook} I. Rubinstein, {\em{Electro-Diffusion of Ions}}, SIAM Studies in Applied Mathematics, SIAM, Philadelphia 1990.
\bibitem{rubinstein}S. M. Rubinstein, G. Manukyan, A. Staicu, I. Rubinstein, B. Zaltzman, R.G.H. Lammertink, F. Mugele, M. Wessling, Direct observation of a nonequilibrium electro-osmotic instability. Phys. Rev. Lett. {\bf{101}}, 
(2008) 236101-236105.
\bibitem{rubizaltz} I. Rubinstein, B. Zaltzman, Electro-osmotically induced convection at a permselective membrane, Phys. Rev. E {\bf{62}} (2000) 2238-2251.
\bibitem {ryham} R. Ryham, Existence, uniqueness, regularity and long-term behavior for dissipative systems modeling electrohydrodynamics. arXiv:0910.4973v1, (2009).
\bibitem{schmuck} M. Schmuck. Analysis of the Navier-Stokes-Nernst-Planck-Poisson system. Math.Models MethodsAppl., {\bf{19}} (2009), 993-1014.
\bibitem{solo1} V.A. Solonnikov, Estimates for solutions of nonstationary Navier-Stokes equations. J. Soviet Math, {\bf{8}} (1977), 213-317.
\bibitem{zaltzrubi}B. Zaltzman, I. Rubinstein, Electro-osmotic slip and electroconvective instability. J. Fluid Mech. {\bf{579}}, (2007) 173-226.
\end{thebibliography}
\end{document}